\documentclass{cocv}

\usepackage[english]{babel}
\usepackage{amsmath, amssymb, graphicx, cite, url}

\def\ds{\displaystyle}
\def\newD{D_1}
\def\newtilde{{\,\,\approx\,\,}}
\def\tmax{t_{\MAX}^1}
\def\tcut{t_{\operatorname{cut}}}
\def\then{\quad\Rightarrow\quad}
\def\tconj{t^1_{\operatorname{conj}}}
\def\cl{\operatorname{cl}}
\newcommand{\vect}[1]{\left( \begin{array}{c} #1 \end{array} \right)}
\newcommand{\lam}{\lambda}
\newcommand{\Diff}{\operatorname{Diff}\nolimits}
\newcommand{\MAX}{\operatorname{MAX}\nolimits}
\newcommand{\Exp}{\operatorname{Exp}\nolimits}
\newcommand{\am}{\operatorname{am}\nolimits}
\newcommand{\sn}{\operatorname{sn}\nolimits}
\newcommand{\cn}{\operatorname{cn}\nolimits}
\newcommand{\dn}{\operatorname{dn}\nolimits}
\newcommand{\E}{\operatorname{E}\nolimits}
\newcommand{\sh}{\operatorname{sinh}\nolimits}
\newcommand{\ch}{\operatorname{cosh}\nolimits}
\newcommand{\tangh}{\operatorname{tanh}\nolimits}
\newcommand{\spann}{\operatorname{span}\nolimits}
\newcommand{\sgn}{\operatorname{sgn}\nolimits}
\def\xSE{\mathop{\rm SE\,}\nolimits}
\def\xSH{\mathop{\rm SH\,}\nolimits}

\begin{document}
\title{Cut time in sub-Riemannian problem on Engel group}
\author{A.~A.~Ardentov}\address{Program Systems Institute of RAS, Pereslavl-Zalessky, Russia; \email{aaa@pereslavl.ru \ yusachkov@gmail.com}}
\author{Yu.~L.~Sachkov}\sameaddress{1}
\date{28.08.2014}
\begin{abstract} The left-invariant sub-Riemannian problem on the Engel group is considered. The problem gives the nilpotent approximation to generic nonholonomic systems in four-dimensional space with two-dimensional control, for instance to a system which describes motion of mobile robot with a trailer. 

The global optimality of extremal trajectories is studied via geometric control theory. 
The global diffeomorphic structure of the exponential mapping is described. As a consequence, 
 the cut time is proved to be equal to the first Maxwell time corresponding to discrete symmetries of the exponential mapping.
\end{abstract}
\begin{resume} ... \end{resume}
\maketitle
\section*{Introduction}
This paper continues the study of the left-invariant sub-Riemannian problem on the Engel group started in~\cite{engel, engel_conj}. This problem is the simplest rank 2 sub-Riemannian problem on a 4-dimensional space: it provides a nilpotent approximation to a generic sub-Riemannian problem of such kind near a generic point.

A sub-Riemannian (SR) structure on a smooth manifold $M$ is a vector distribution

$$\Delta = \{\Delta_q \subset T_q M ~ | ~ q \in M\} \subset TM$$
with a scalar product in $\Delta$:

$$g = \{g_q \text{ --- scalar product in } \Delta_q ~ | ~ q \in M \}.$$
The subspaces $\Delta_q \subset T_q M$ and the scalar product $g_q\colon \Delta_q \times \Delta_q \to \mathbb{R}$
depend smoothly on a point $q \in M$. The dimension of the subspaces $\Delta_q$ is constant ($\dim \Delta_q$ is called the rank of the distribution $\Delta$).

A Lipschitz curve $q\colon [0, t_1] \to M$ is horizontal if $\dot{q}(t) \in \Delta_{q(t)}$ for almost all $t \in [0, t_1]$. The length of a horizontal curve is 
$l = \int\limits_0^{t_1} g\big(\dot{q}(t), \dot{q}(t)\big)^{1/2} dt$. The sub-Riemannian distance between points $q_0$, $q_1 \in M$ is the infimum of lengths of horizontal curves that connect $q_0$ to $q_1$.
A horizontal curve $q(t), t \in [0, t_1]$, is a (length) minimizer if it has a minimum possible length among all horizontal curves that connect the points $q(0)$ and $q(t_1)$. Description of minimizers is one of important problems of sub-Riemannian geometry. The most efficient approach to this problem is given by geometric control theory~\cite{jurd, notes, abb}, it consists of the following steps:
\begin{enumerate}
\item\label{step1} proof of existence of minimizers,
\item\label{step2} description of SR geodesics (\ie, curves whose small arcs are minimizers),
\item\label{step3} selection of minimizers among geodesics.
\end{enumerate}

Step~(\ref{step1}) is straightforward. If $M$ is connected and $\Delta$ is bracket generating, \ie, $\mathrm{Lie}_q \Delta = T_q M, \quad \forall q \in M,$
then any points $q_0, q_1 \in M$ can be connected one to another by a horizontal curve (Rashevsky-Chow theorem). If additionally the point $q_1$ is sufficiently close to $q_0$, or if the SR distance is complete, or if $\Delta$ and $g$ are left-invariant on a Lie group $M$, then $q_0$ can be connected with $q_1$ by a minimizer (Filippov theorem).

Step~(\ref{step2}) is performed via application of Pontryagin maximum principle (PMP), which states that any geodesic (thus any minimizer) is a projection of a trajectory of a certain Hamiltonian system on the cotangent bundle $T^* M$. So the second step reduces to the study of integrability of the Hamiltonian system of PMP and efficient parameterization of trajectories of this system.

Step~(\ref{step3}) is the hardest one. Local optimality of geodesics (\ie, optimality w.r.t. sufficiently close geodesics) is studied via conjugate points estimates. For the study of global optimality in problems with a big symmetry group, one can often obtain bounds (or explicit description) of cut time via the study of symmetries and global structure of the exponential mapping. We suggest the following detailing of Step~(\ref{step3}) first applied in~\cite{abck} and further developed in~\cite{dido_exp, max1, max2, max3, max_sre, cut_sre1, cut_sre2, el_max, el_conj, el_exp}:
\begin{itemize}
\item[(3.1)] Discrete and continuous symmetries of the exponential mapping are found;
\item[(3.2)] Maxwell points corresponding to the symmetries are found (\ie, points where several geodesics obtained one from another by a symmetry meet one another). These points (and their preimage via exponential mapping) form the Maxwell strata in the image (resp., in the preimage) of the exponential mapping. Along each geodesic, the first Maxwell time corresponding to the symmetries (\ie, the time when the geodesic meets a Maxwell strata) is found;
\item[(3.3)] One proves that for any geodesic the first conjugate time is greater or equal to the first Maxwell time corresponding to the symmetries. Here the homotopy invariance of Maslov index (number of conjugate points on a geodesic) can be applied~\cite{cime};
\item[(3.4)] One considers restriction of the exponential mapping to the subdomains cut out in preimage and image of this mapping by the Maxwell strata corresponding to symmetries, and proves that this restriction is a diffeomorphism via Hadamard global diffeomorphism theorem~\cite{Hadamard}; 
\item[(3.5)] On the basis of the global structure of the exponential mapping thus described, it is often possible to prove that the cut time along a geodesic (\ie, time when it loses its global optimality) is equal to the first Maxwell time corresponding to symmetries. Moreover, in this way one proves that for any terminal point in a subdomain in the image of the exponential mapping, there exists a unique minimizer which can be computed by inverting the exponential mapping in the subdomain;
\item[(3.6)] Finally, for systems with big symmetry group one can construct the full optimal synthesis, and numerical algorithms and software for computation of optimal trajectories with given boundary conditions.
\end{itemize}

So far, the approach described has been applied in full just to several problems: SR problem in the flat Martinet case~\cite{abck},
SR problems on $\xSO(3)$ and $\xSL(2)$ with the Killing metric~\cite{boscainSO3},
SR problem on $\xSE(2)$~\cite{max_sre, cut_sre1, cut_sre2},
Euler elastic problem~\cite{el_max, el_conj, el_exp}.
There are partial results on the nilpotent SR problem with the growth vector (2,~3,~5)~\cite{dido_exp, max1, max2, max3}
and SR problem on $\xSH(2)$~\cite{sh2_exp,sh2_conj}.

For the SR problem on the Engel group, Step~(\ref{step1}), Step~(\ref{step2}) and Steps (3.1), (3.2) are performed in~\cite{engel}
while (3.3) is done in~\cite{engel_conj}.
The aim of this paper is to perform Steps (3.4), (3.5). We recall the results previously obtained in the next section.

The sub-Riemannian problem on the Engel group is a left-invariant problem on a Lie group. Such problems receive significant attention in geometric control since they provide very symmetric models which can often be studied explicitly in great detail. For left-invariant SR problems on Lie groups, one can often describe optimal synthesis, the structure of spheres, cut and conjugate loci. This information can give insight for general problems, where such a detailed study is much more complicated.

Left-invariant SR problems on 3D and 4D Lie groups have recently been fully classified \cite{agrachev_barilari, almeida}.
In the 3-dimensional case, optimal synthesis is known for the Heisenberg group~\cite{vershik-gershkovich},
for $\xSO(3)$ and $\xSL(2)$ with the Killing metric~\cite{boscainSO3} and for $\xSE(2)$~\cite{max_sre, cut_sre1, cut_sre2}.
This work continues a detailed study of the simplest 4-dimensional case.

\section{Previously obtained results}
In this section we recall results on the SR problem on the Engel group obtained previously in works~\cite{engel, engel_conj}.
\subsection {Problem statement}
The Engel group is the 4-dimensional Lie group represented by matrices as follows:

\begin{align*}
&M = \left\{
\left(
\begin{array}{cccc}
1 & b & c & d \\
0 & 1 & a & a^2/2 \\
0 & 0 & 1 & a \\
0 & 0 & 0 & 1 
\end{array}
\right)
\mid
a, \ b, \ c, \ d \in \xR
\right\}.
\end{align*}
It is a 4-dimensional nilpotent Lie group, connected and simply connected (see an explanation of the name "Engel" for this group in~\cite{mont}, Sec.~6.11).

The Lie algebra of the Engel group is the four-dimensional nilpotent Lie algebra $L = \spann(X_1, X_2, X_3, X_4)$ with the multiplication table

\begin{minipage}[h]{0.2\linewidth}
\centering \includegraphics[width=0.75\linewidth]{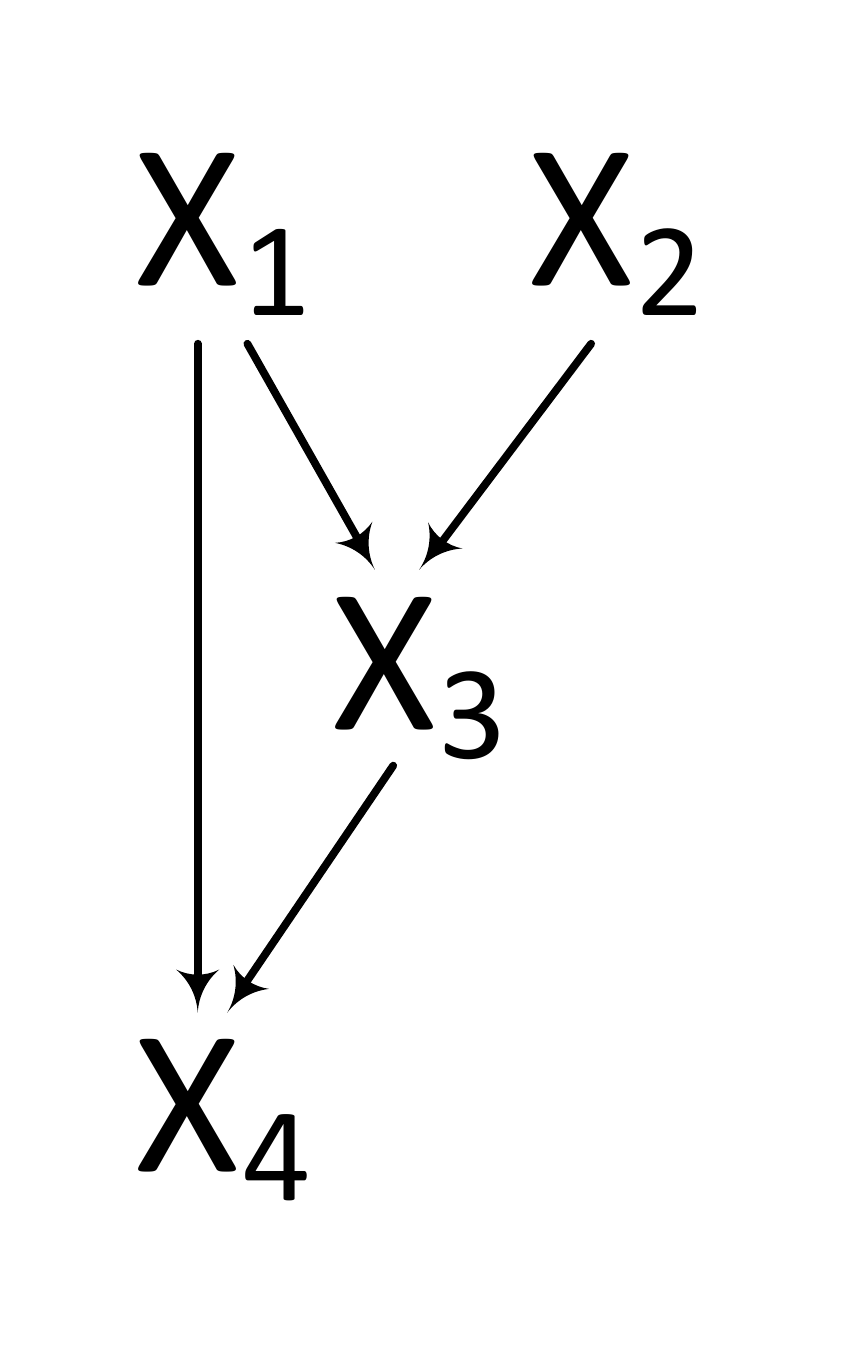}
\end{minipage}
\begin{minipage}[h]{0.772\linewidth}
\begin{center}
\begin{align}
&X_3=[X_1, X_2], \nonumber \\
&X_4=[X_1, X_3], \nonumber\\ 
&[X_2,X_3]=[X_1,X_4]=[X_2,X_4]=0. \label{mult}
\end{align}

~
\end{center}
\end{minipage} \\
Thus it has graduation
\begin{align*}
& L = L_1 \oplus L_2 \oplus L_3,\\
& L_1 = \spann(X_1, X_2), \quad L_2 = \mathbb{R} X_3, \quad L_3 = \mathbb{R} X_4, \\
&[L_i, L_j] = L_{i+j}, \ L_k = \{0\} \text{ for } k  \geq 4,
\end{align*}
and the Engel group is a Carnot group~\cite{mont}.

We consider the sub-Riemannian problem on the Engel group $M$ for the left-invariant sub-Riemannian structure generated by the orthonormal frame $X_1$, $X_2$:
\begin{align*}
&\dot q = u_1 X_1(q) + u_2 X_2(q), \qquad q \in M, \quad (u_1, u_2) \in \xR^2, \\
& q(0) = q_0, \quad q(t_1) = q_1, \\
& l = \int_0^{t_1} \sqrt{u_1^2 + u_2^2} \, d t \rightarrow \min.
\end{align*}

In appropriate coordinates $q = (x,y,z,v)$ on the Engel group $M \cong \xR^4$, the problem 
 is stated as follows:
\begin{align} 
&\dot{q} = \vect{\dot{x} \\ \dot{y} \\ \dot{z} \\ \dot{v}} = u_1 
\vect{1 \\ 0 \\ - y/2 \\ 0} + u_2 \vect{0 \\ 1 \\ x/2 
\\ (x^2+y^2)/2}, \quad q = (x,y,z,v) \in M=\xR^4, \quad (u_1, u_2) \in \xR^2, \label{pr1}\\
& q(0) = q_0 = (x_0,y_0,z_0,v_0), \quad q(t_1) = q_1 = (x_1, y_1, z_1, v_1), \label{pr2} \\
& l = \int_0^{t_1} \sqrt{u_1^2 + u_2^2} \, d t \rightarrow \min.
\label{pr3}
\end{align}

By virtue of the multiplication table~(\ref{mult}) for the vector fields of the orthonormal frame
$$ X_1 = \frac{\partial}{\partial x} - \frac{y}{2} \frac{\partial}{\partial z}, \qquad  X_2 = \frac{\partial}{\partial y} +\frac {x}{2} \frac{\partial}{\partial z} + \frac{x^2 + y^2}{2} \frac{\partial}{\partial v}$$
and their Lie brackets
$$ X_3 = [X_1, X_2] = \frac{\partial}{\partial z} + x \frac{\partial}{\partial v} , \qquad X_4 = [X_1, X_3] = \frac{\partial}{\partial v}, $$
system~(\ref{pr1}) is completely controllable, \ie, any points $q_0, q_1 \in \mathbb{R}^4$ can be connected by its trajectory.

Since the problem is invariant under left shifts on the Engel group, we can assume that the initial point is the identity $q_0 = (x_0, y_0, z_0, v_0) = (0, 0, 0, 0)$.

\subsection {Parameterization of geodesics}
Existence of optimal solutions of problem~(\ref{pr1})--(\ref{pr3}) is implied by Filippov theorem~\cite{notes}. By the CauchyЦ-Schwarz inequality, it follows that sub-Riemannian length minimization problem~(\ref{pr3}) is equivalent to the action minimization problem:
\begin{align} \label{J}
\int_0^{t_1} \frac{u_1^2+u_2^2}{2} \, d t \rightarrow \min.
\end{align}
Pontryagin maximum principle~\cite{PGBM, notes} was applied to the resulting optimal control problem~(\ref{pr1}), (\ref{pr2}), (\ref{J}) in~\cite{engel}. 

A sub-Riemannian geodesic can be normal or abnormal, or both. For the SR problem on the Engel group, each abnormal geodesic is simultaneously normal, thus in the sequel we consider only normal geodesics. 

Normal geodesics are projections $q_t = \pi ({\lambda}_t)$ via the canonical projection $\pi \colon T^* M \to M$ of solutions to the Hamiltonian system
\begin{align}\label{norm_ham}
 \dot{\lambda} = \vec{H} (\lambda), \quad \lambda \in T^* M, 
\end{align}
with the Hamiltonian function $H = \frac{1}{2} (h^2_1 + h^2_2).$ Here and below $h_i (\lambda) = \left\langle \lambda, X_i (q) \right\rangle, \ \lambda \in T^* M, \ i = 1, \dots , 4,$ are Hamiltonians that correspond to the left-invariant frame and are linear on fibers of the cotangent bundle $T^* M$.

Arclength parameterized geodesics (\ie, with velocity $g(\dot{q_t},\dot{q_t}) \equiv 1$) are projections of extremals ${\lambda}_t$ lying on the level surface $\{ \lambda \in T^* M ~ | ~ H(\lambda) = 1/2 \}$.

\begin{figure}[htbp]
\centering
\includegraphics[width=0.47\linewidth]{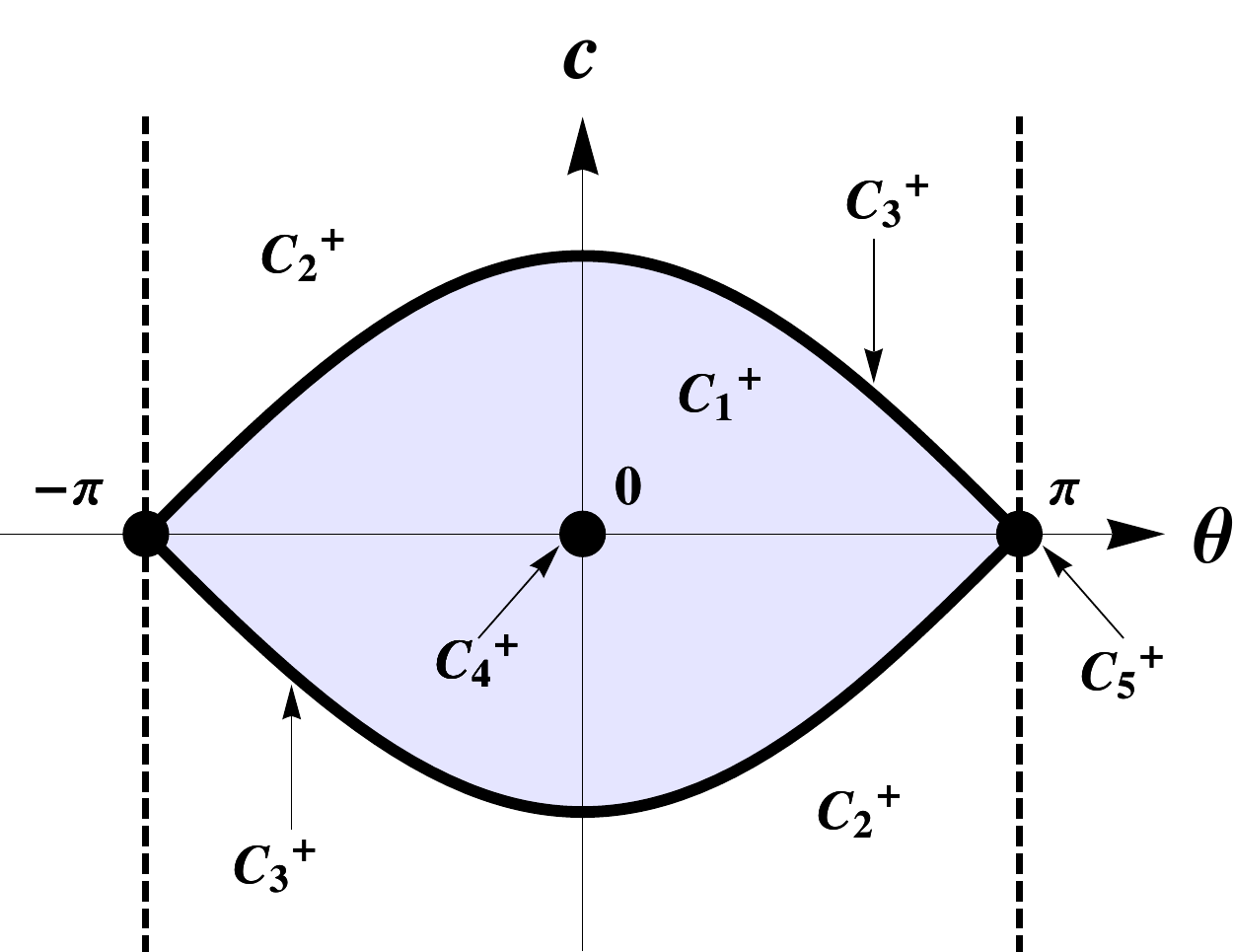}
\qquad
\includegraphics[width=0.47\linewidth]{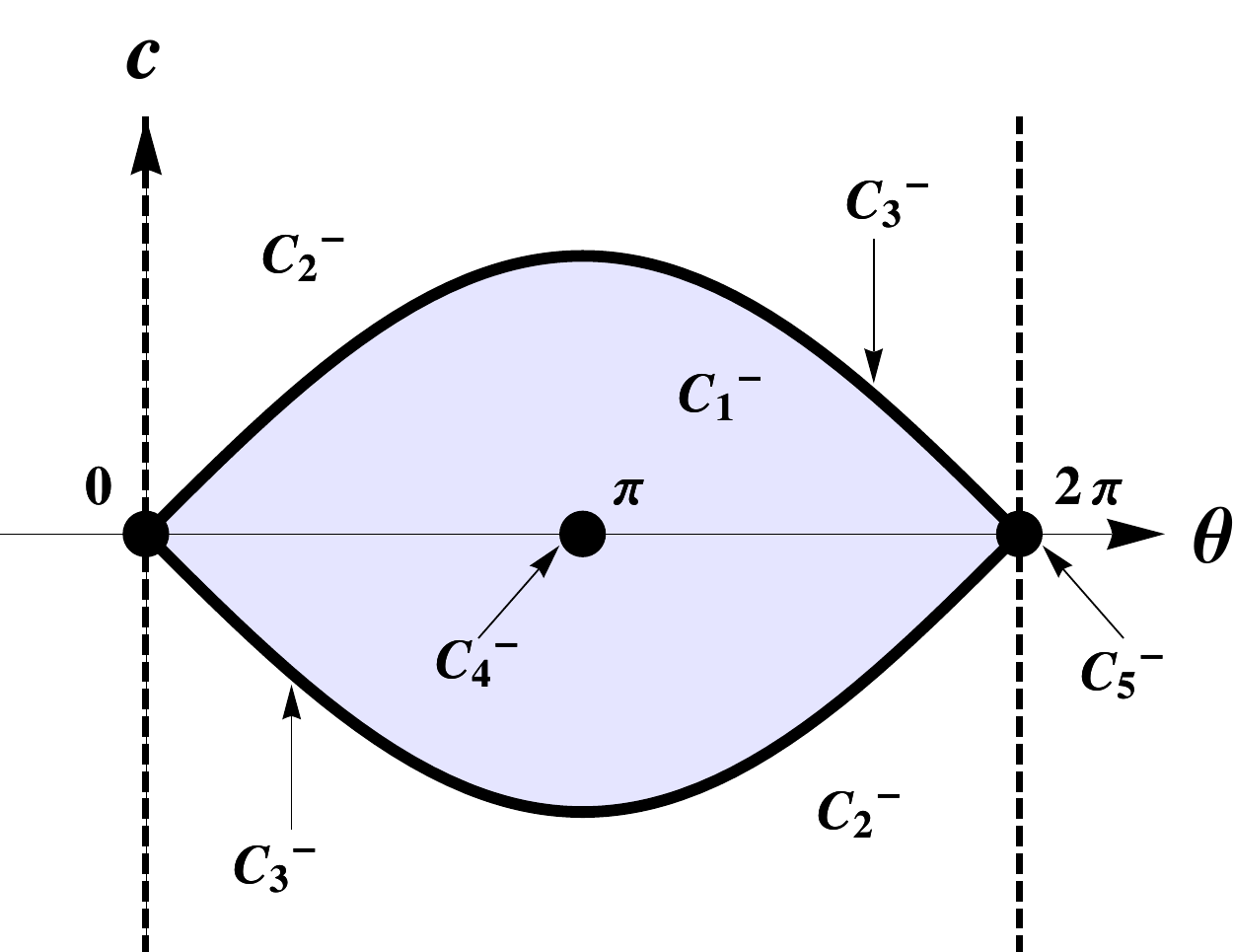}
\caption{Stratification of $C$ for $\alpha>0$ and for $\alpha<0$}\label{pendulum+-}
\end{figure}

Introduce coordinates $(\theta, c, \alpha)$ on the level surface $\left\{\lam \in T^* M \mid H=1/2\right\}$ by the following formulas:
$$
h_1 = \cos (\theta+  \pi/2), \qquad h_2 = \sin (\theta+ \pi/2), \qquad
h_3 = c, \qquad h_4 = \alpha.
$$
On this surface
the normal Hamiltonian system~(\ref{norm_ham}) takes the following form:
\begin{eqnarray}
&&\dot{\theta}= c, \qquad \dot{c}= - \alpha\, \sin \theta, \qquad \dot{\alpha}=0, \label{pend}\\ 
&&\dot{q}= \cos \theta \, X_1(q)+\sin\theta \, X_2 (q), \qquad q(0)=q_0. \nonumber
\end{eqnarray}
The family of all normal extremals is parameterized by points of the phase cylinder of pendulum
\begin{eqnarray*}
C = \left\{\lambda \in T_{q_0} ^* M \mid H(\lambda)= 1/2\right\}= \left\{(\theta, c, \alpha  ) \mid \theta \in S^1, \ c, \alpha \in \xR  \right\},
\end{eqnarray*}
and is given by the exponential mapping
\begin{eqnarray*}
&&\Exp \colon N = C \times \xR_+ \to M,\\
&&\Exp (\lambda, t) = q_t = (x_t, y_t, z_t, v_t).
\end{eqnarray*}
The energy integral of pendulum~(\ref{pend}) is given by $\ds E=\frac{c^2}{2}-\alpha \cos \theta$. The cylinder $C$ has the following stratification corresponding to the particular type of trajectories of the pendulum:
\begin{align}
&C=\cup_{i=1}^7 C_i,   \quad  C_i \cap C_j = \emptyset, \ i \neq j, 
\quad \lambda = (\theta,c,\alpha), \nonumber \\
&C_1 = \{\lambda \in C \mid \alpha \neq 0, E\in(- |\alpha|, 
|\alpha|)\}, \label{C1}\\
&C_2 = \{\lambda \in C \mid \alpha \neq 0, E\in(|\alpha|,+\infty)\}, \label{C2}
\\
&C_3 = \{\lambda \in C \mid \alpha \neq 0, E=|\alpha|, c \neq 0 \}, \label{C3}\\
&C_4 = \{\lambda \in C \mid \alpha \neq 0, E=-|\alpha|\}, \label{C4}\\
&C_5 = \{\lambda \in C \mid \alpha \neq 0, E=|\alpha|, c = 0\}, \label{C5}\\
&C_{6} = \{\lambda \in C \mid \alpha = 0, \ c \neq 0\}, \label{C6}\\
&C_7 = \{\lambda \in C \mid \alpha = c = 0\}. \label{C7}
\end{align}
Further, the sets $C_i, \, i=1, \dots, 5,$ are divided into subsets determined by the sign of $\alpha$ (see Fig.~\ref{pendulum+-}):
\begin{align*}
&C_i^+ = C_i \cap \{\alpha>0\}, \qquad C_i^- = C_i \cap 
\{\alpha<0\}, \qquad i = 1, \dots, 5.
\end{align*}

In order to parameterize extremal trajectories, coordinates $(\varphi, k, \alpha)$ in the domains $C_1$ and $C_2$ were introduced in~\cite{engel} in the following way.
 
In the domain $C_1^+$
\begin{align*}
&k=\sqrt{\frac{E+\alpha}{2 \alpha}} = \sqrt{\frac{c^2}{4 
\alpha}+\sin^2 \frac{\theta}{2}}\in (0,1), \\
&\sin\frac{\theta}{2} = k \sn (\sqrt{\alpha} \varphi), \quad
\cos\frac{\theta}{2} = \dn (\sqrt{\alpha} \varphi), \quad
\frac {c}{2} = k \sqrt{\alpha} \cn (\sqrt{\alpha} \varphi), \quad \varphi 
\in [0,4 K(k)].
\end{align*}

In the domain $C_2^+$
\begin{align*}
&k=\sqrt{\frac{2 \alpha}{E+\alpha}} = \frac{1}{\sqrt{\frac{c^2} {4 
\alpha}+\sin^2 \frac{\theta}{2} }}\in (0,1),\\
&\sin\frac{\theta}{2} = \sgn c \, \sn \frac{\sqrt{\alpha} \varphi}{k}, \quad
\cos\frac{\theta}{2} = \cn \frac{\sqrt{\alpha} \varphi}{k}, \quad
\frac{c}{2} = \sgn c \, \frac {\sqrt{\alpha}}{k}  \dn \frac 
{\sqrt{\alpha} \varphi}{k}, \quad \varphi \in [0,2 k K(k)]. 
\end{align*}
Here and below $\dn$, $\sn$, $\cn$ are Jacobi elliptic functions depending on $k$, $K(k)$ is the complete elliptic integral of the first kind~\cite{whit_vatson}.

In the domains $C_1^-$, $C_2^-$ the coordinates $\varphi$ and $k$ are defined as follows:  
$$
\varphi (\theta, c, \alpha) = \varphi (\theta - \pi, c, -\alpha), 
\qquad
k(\theta, c, \alpha) = k (\theta-\pi, c, -\alpha). 
$$
Immediate differentiation shows that system~(\ref{pend}) rectifies in the coordinates $(\varphi, k, \alpha)$: $$\dot{\varphi}=1, \qquad \dot{k}=0, \qquad \dot{\alpha}=0.$$

In terms of these coordinates, geodesics $q_t = \Exp(\lambda , t)$ with $\lambda = (\theta , c, \alpha) \in \bigcup_{i=1}^{3}C_i$ and $\alpha = 1$ are parameterized as follows.

If $\lambda \in C_1$, then
\begin{align}
&x_t = 2 k (\cn \varphi_t - \cn \varphi), \nonumber \\
&y_t = 2 \big(\E(\varphi_t) - \E(\varphi)\big)-t, \nonumber\\
&z_t = 2 k \big(\sn \varphi_t \dn \varphi_t - \sn \varphi \dn \varphi 
- \frac{y_t}{2} (\cn \varphi_t + \cn \varphi)\big), \nonumber\\
&v_t= \frac{y_t^3}{6} + 2 k^2 \cn^2 \varphi y_t - 4 k^2 \cn \varphi 
(\sn \varphi_t \dn \varphi_t - \sn \varphi \dn \varphi)  \nonumber \\
&\quad+ 2 k^2 \bigg(\frac{2}{3} \cn \varphi_t \dn \varphi_t \sn 
\varphi_t- \frac{2}{3} \cn \varphi \dn \varphi \sn \varphi + \frac{1-
k^2}{3 k^2} t  +\frac{2 k^2 -1}{3 k^2}\Big(\E(\varphi_t)-\E(\varphi)\Big)\bigg). \label{ExpC1}
\end{align}
Here and below $\E(\varphi) = \int_0^\varphi \dn^2 t \xdif t = E \big(\am (\varphi), k\big)$ is the Jacobi epsilon function and $E(u,k)$ is incomplete elliptic integral of the second kind depending on Jacobi amplitude $u = \am(\varphi)$ and $k$. The Jacobi amplitude $\am (\varphi)$ is the inverse function of the incomplete elliptic integral of the first kind: $F\big(\am(\varphi)\big) = \varphi$. 

\medskip\noindent
If $\lambda \in C_2$, then
\begin{align}
&x_t = \frac {2 \sgn c} {k} \Big(\dn \psi_t - \dn \psi \Big), \nonumber\\
&y_t = \frac {k^2-2}{k^2} t + \frac {2}{k} \Big(\E(\psi_t) - \E(\psi)\Big), \nonumber\\
&z_t = - \frac{x_t y_t}{2} - \frac{2 \sgn c \dn \psi}{k} y_t + 2 \sgn c \, (\cn \psi_t \sn \psi_t - \cn \psi \sn \psi), \nonumber\\
&v_t= \frac {4}{k} \bigg( \frac {1}{3}\cn \psi_t \dn \psi_t \sn \psi_t 
- \frac{1}{3} \cn \psi \dn \psi \sn \psi -\frac{1-k^2}{3 k^3} t - \frac {k^2-2}{6 k^2} \Big(\E(\psi_t)-\E(\psi)\Big)\bigg) \nonumber\\
&\quad + \frac {y_t^3}{6} + \frac {2 y_t}{k^2} \dn^2 \psi - \frac{4}{k} \dn \psi \big(\cn \psi_t \sn \psi_t - \cn \psi \sn 
\psi \big), \nonumber \\
&\psi=\frac{\varphi}{k}, \quad \psi_t = \psi + \frac{t}{k}. \label{ExpC2}
\end{align}

\medskip\noindent
If $\lambda \in C_3$, then
\begin{align}
&x_t = 2 \sgn c \left(\frac {1}{\ch \varphi_t} - \frac{1}{\ch \varphi}\right), \nonumber\\
&y_t = 2 (\tangh \varphi_t - \tangh \varphi) - t , \nonumber\\
&z_t = - \frac{x_t y_t}{2} - \frac{2 \sgn c}{\ch \varphi} y_t + 2 \sgn 
c\left(\frac{\tangh \varphi_t}{\ch \varphi_t} - \frac{\tangh \varphi}{\ch 
\varphi}\right), \nonumber \\
&v_t= \frac{2}{3}\left(\tangh \varphi_t - \tangh \varphi + 2 \frac{\tangh 
\varphi_t}{\ch^2 \varphi_t}-2\frac{\tangh \varphi}{\ch^2 \varphi}\right) 
+ \frac{y_t^3}{6} + \frac {2 y_t}{\ch^2 \varphi}  -  \frac{4}{\ch \varphi} \left(\frac{\tangh \varphi_t}{\ch 
\varphi_t} - \frac{\tangh \varphi}{\ch \varphi}\right). \label{ExpC3}
\end{align}

Parameterization of geodesics for $\lambda \in \bigcup_{i=1}^{3}C_i$ and arbitrary $\alpha \ne 0$ is obtained from the above parameterization for $\alpha = 1$ via the following symmetries of the Hamiltonian system: dilations
\begin{align*}
& \delta_\mu : (\theta, c, \alpha, t, x, y, z, v) ~ \mapsto ~ (\theta, c/\mu, \alpha/\mu^2, \mu t, \mu x, \mu y, \mu^2 z, \mu^3 v),\quad \mu > 0, \\
& \delta_\mu : (\varphi, k, \alpha) ~ \mapsto ~ (\mu \varphi, k, \alpha/\mu^2), 
\end{align*}
and reflection
\begin{align*}
&(\theta, c, \alpha, t, x, y, z, v) ~ \mapsto ~ (\theta - \pi, c, -\alpha, t, -x, -y, z, -v), \\
&(\varphi, k, \alpha) ~ \mapsto ~ (\varphi, k, -\alpha). 
\end{align*}

In the remaining cases $\lambda \in \bigcup_{i=4}^{7}C_i$ geodesics are parameterized by elementary functions as follows.

$\lambda \in C_4:$
\begin{align}
&x_t = 0, \qquad y_t = t \sgn \alpha , \qquad z_t = 0, \qquad v_t = 
\frac{t^3}{6} \sgn \alpha. \label{ExpC4}
\end{align}

$\lambda \in C_5:$
\begin{align}
&x_t = 0, \qquad y_t = - t \sgn \alpha, \qquad z_t = 0, \qquad v_t = - 
\frac{t^3}{6} \sgn \alpha. \label{ExpC5}
\end{align}

$\lambda \in C_6:$
\begin{align}
&x_t = \frac{\cos (c t + \theta) - \cos \theta}{c}, &&y_t = \frac{\sin 
(c t + \theta) - \sin \theta}{c}, \nonumber\\
&z_t = \frac{c t - \sin (c t)}{2 c^2}, 
&&v_t = \frac{3 \cos \theta - 2 c t \sin \theta - 4 \cos (c t + \theta) + \cos (2 c t + \theta)}{4 c^3}.  \label{ExpC6}
\end{align}

$\lambda \in C_7:$
\begin{align}
&x_t = - t \sin \theta, \qquad y_t = t \cos \theta, \qquad z_t = 0, 
\qquad v_t = \frac{t^3}{6} \cos \theta. \label{ExpC7}
\end{align}

Projections of geodesics to the plane $(x, y)$ are Euler elasticae (stationary configurations of planar elastic rod with fixed endpoints and tangents at endpoints)~\cite{euler, love, el_max, el_conj, el_exp}: inflexional ones for $\lambda \in C_1$, non-inflexional ones for $\lambda \in C_2$, critical ones for $\lambda \in C_3$, straight lines for $\lambda \in C_4 \cup C_5 \cup C_7$, and circles for $\lambda \in C_6$.

\subsection{Symmetries of exponential mapping}

A pair of mappings 
$$ s \colon N \to N, \qquad s \colon M \to M $$
is called a symmetry of the exponential mapping if it commutes with this mapping:
$$ s \circ \Exp(\lambda, t) = \Exp \circ s(\lambda, t), \quad (\lambda, t) \in N. $$

\subsection{Dilations}

A one-parameter group of symmetries of the exponential mapping is formed by dilations
\begin{align}
& \delta_\mu : (\theta, c, \alpha, t) ~ \mapsto ~ (\theta, c/\mu, \alpha/\mu^2, \mu t), \nonumber \\
& \delta_\mu : (x, y, z, v) ~ \mapsto ~ (\mu x, \mu y, \mu^2 z, \mu^3 z), \quad \mu > 0. \label{dmuM}
\end{align}

\subsection{Reflections}

The following mappings ${\varepsilon}^i \colon C \to C$ preserve the field of directions of the vertical part of the Hamiltonian vector field $\vec{H}_v = c \frac{\partial}{\partial \theta} - \alpha \sin \theta \frac{\partial}{\partial c} \in \mathrm{Vec}(C)$:
\begin{align*}
&\varepsilon^1 : (\theta,c, \alpha) ~ \mapsto ~ (\theta, -c, \alpha), &&\varepsilon^2 : (\theta,c, \alpha) ~ \mapsto ~ (-\theta, c, \alpha), \\
&\varepsilon^3 : (\theta,c, \alpha) ~ \mapsto ~ (-\theta, -c, \alpha),  &&\varepsilon^4 : (\theta,c, \alpha) ~ \mapsto ~ (\theta+\pi, c, -\alpha), \\
&\varepsilon^5 : (\theta,c, \alpha) ~ \mapsto ~ (\theta+\pi, -c, -\alpha), &&\varepsilon^6 : (\theta,c, \alpha) ~ \mapsto ~ (-\theta+\pi, c, -\alpha), \\
&\varepsilon^7 : (\theta,c, \alpha) ~ \mapsto ~ (-\theta+\pi, -c, -\alpha). 
\end{align*}
More precisely, ${\varepsilon}^i_* \vec{H}_v = \vec{H}_v$ for $i = 3, 4, 7,$ and ${\varepsilon}^i_* \vec{H}_v = - \vec{H}_v$ for $i = 1, 2, 5, 6$. The action of reflections ${\varepsilon}^i$ is continued to symmetries of the exponential mapping as follows. 

The action ${\varepsilon}^i \colon N \to N$ is defined as

\[
{\varepsilon}^i (\lambda, t) = 
 \begin{cases}
 \big({\varepsilon}^i (\lambda), t\big), & \textrm{if }{\varepsilon}^i_*  \vec{H}_v = \vec{H}_v, \\
 \big({\varepsilon}^i \circ e^{t \vec{H}_v} (\lambda), t\big), & \textrm{if } {\varepsilon}^i_* \vec{H}_v = - \vec{H}_v.
\end{cases}
\]

The action ${\varepsilon}^i \colon M \to M$ is defined as
\begin{align}
&\varepsilon^i(q) = \varepsilon^i(x,y,z,v) = q^i = (x^i, y^i, z^i, v^i), \label{epsiM}\\
&(x^1, y^1, z^1, v^1) = (x, y, -z, v - x z), \label{eps1M}\\
&(x^2, y^2, z^2, v^2) = (-x, y, z, v - x z), \label{eps2M}\\
&(x^3, y^3, z^3, v^3) = (-x, y, -z, v), \label{eps3M}\\
&(x^4, y^4, z^4, v^4) = (-x, -y, z, -v), \label{eps4M}\\
&(x^5, y^5, z^5, v^5) = (-x, -y, -z, -v + x z), \label{eps5M}\\
&(x^6, y^6, z^6, v^6) = (x, -y, z, -v + x z), \label{eps6M}\\
&(x^7, y^7, z^7, v^7) = (x, -y, -z, -v). \label{eps7M}
\end{align}

Thus defined reflections ${\varepsilon}^i , i = 1, \ldots, 7,$ form a discrete group of symmetries of the exponential mapping (together with the identity mapping).

\subsection{Maxwell points}
A point $q_t$ of an extremal trajectory $q_s = \Exp(\lam,s)$ is called a Maxwell point if there exists another extremal trajectory $\widetilde{q}_s = \Exp(\widetilde \lam, s)$, $\widetilde{q}_s \not\equiv q_s$, such that $\widetilde{q}_t = q_t$. The instant $t$ is called a Maxwell time. It is known~\cite{max3} that an extremal trajectory cannot be optimal after a Maxwell time.

The main result of paper~\cite{engel}, given by Th.~\ref{th:tcut_bound} below, provides an upper bound of the cut time along extremal curves
$$\tcut(\lambda) = \sup \{ t>0 \mid \mathrm {Exp} (\lambda, s) \text{ is optimal for } s \in [0,t]\}.$$
Define the following function $\tmax \colon C \to (0, +\infty]$:
\begin{align}
&\lambda \in C_1 & \then &\tmax = \min \big(2 p_z^1(k), 4 K(k)\big)/\sigma, \label{tmaxC1}\\
&\lambda \in C_2 & \then &\tmax = 2 k K(k) /\sigma, \label{tmaxC2}\\
&\lambda \in C_6 & \then &\tmax = 2 \pi / |c|, \label{tmaxC6}\\
&\lambda \in C_3 \cup C_4 \cup C_5 \cup C_7 & \then &\tmax=+\infty. \label{tmaxC4}
\end{align}
where $\sigma = \sqrt{|\alpha|}$; $\ds K(k)= \int_0^\frac{\pi}{2} \frac{dt}{\sqrt{1- k^2 \sin^2 t}};$
$p^1_z(k)\in \big(K(k), 3K(k)\big)$ is the first positive root of the function $f_z(p,k)=\dn p \,\sn p+ (p-2\E(p))\cn p;$ 

\begin{thrm}[\cite{engel}, Th. 3]\label{th:tcut_bound}
For any $\lambda \in C$ 
\begin{align}
\tcut(\lambda) \leq \tmax(\lambda). \label{tcutbound}
\end{align}
\end{thrm}

\begin{prpstn} \label{tMaxinvar}
The function $\tmax \colon C \to (0, +\infty]$ has the following invariant properties:
\begin{enumerate}
\item\label{inv1} $\tmax (\lambda)$ depends only on the values of $E$ and $|\alpha|$,
\item\label{inv2} $\tmax (\lambda)$ is an integral of the vector field $\vec{H}_v$,
\item\label{inv3} $\tmax (\lambda)$ is invariant w.r.t. reflections: if $\lambda \in C, \ {\lambda}^i  = {\varepsilon}^i (\lambda) \in C$, then $\tmax ({\lambda}^i) = \tmax (\lambda)$,
\item\label{inv4} $\tmax$ respects the action of dilations: if $\lambda \in C, \ {\lambda}_\mu = \delta_\mu (\lambda),$ then $\tmax ({\lambda}_\mu) = \mu \tmax (\lambda)$.
\end{enumerate}
\end{prpstn}

\begin{proof}
(\ref{inv1}) We denote by $\sqcup$ the union of disjoint sets. Notice first that the decomposition 
\begin{align} \label{C1357}
C = C_1 \sqcup C_2 \sqcup C_{35} \sqcup C_4 \sqcup C_6 \sqcup C_7 
\end{align}
with $C_{35} = C_3 \cup C_5 = \{ \lambda \in C ~ | ~ \alpha \ne 0, E = |\alpha| \}$ is determined only by the functions $E$ and $|\alpha|$, see definitions~(\ref{C1})--(\ref{C7}).
Thus it remains to show that restriction
of $\tmax$ to each of the subsets in decomposition~(\ref{C1357}) depends only on $E$ and $|\alpha|$.

If $\lambda \in C_1$, then $\ds k = \sqrt{\frac{E + |\alpha|}{2 |\alpha|}}$, thus $k = k(E, |\alpha|)$, so $\tmax = \tmax (E, |\alpha|)$. 

The case $\lambda \in C_2$ is similar to the case $\lambda \in C_1$. 

If $\lambda \in C_{35} \cup C_4 \cup C_7$, then $\tmax = +\infty = \tmax (E, |\alpha|)$.

Finally, if $\lambda \in C_6$, then $\ds \tmax = \frac{2 \pi}{ |c| } = \frac{\sqrt{2} \pi}{ \sqrt{E} }$.

\medskip (\ref{inv2}) Since $E$ and $\alpha$ are integrals of the vector field $\vec{H}_v$, then $\tmax = \tmax (E, |\alpha|)$ is an integral of $\vec{H}_v$ as well.

\medskip (\ref{inv3}) Let $\lambda \in C, \ {\varepsilon}^i (\lambda) = {\lambda}^i \in C$. Since $E({\lambda}^i) = E(\lambda), \ \alpha({\lambda}^i) = \pm \alpha (\lambda)$ and $\tmax = \tmax (E, |\alpha|)$, then $\tmax ({\lambda}^i) = \tmax (\lambda)$.

\medskip (\ref{inv4}) Let $\lambda \in C $ and ${\lambda}_\mu = \delta_\mu (\lambda), \ \mu > 0$. Since we have $\ds E ({\lambda}_\mu) = \frac{1}{\mu^2} E(\lambda)$ and $\ds \alpha ({\lambda}_\mu) = \frac{1}{\mu^2} \alpha (\lambda),$
then $\delta_\mu (C_i) = C_i, \ i = 1, \dots, 7$, and $k({\lambda}_\mu) = k(\lambda)$. Then it follows from the definition of the function $\tmax$ that $\tmax ({\lambda}_\mu) = \mu \cdot  \tmax (\lambda)$  for $\lambda \in C_i$ and each $i = 1, \dots, 7$.
\end{proof}

\subsection{Conjugate points}
A point $q_t = \Exp(\lambda, t) $ is called a conjugate point for $q_0$ if $\nu = (\lambda, t)$ is a critical point of the exponential mapping and that is why $q_t$ is the corresponding critical value:
$$
d_{\nu} \Exp \colon T_{\nu}N \to T_{q_t}M \text{ is degenerate}.
$$
The instant $t$ is called a conjugate time along the extremal trajectory $q_s = \Exp(\lambda, s)$, $s\geq 0$.

The first conjugate time along a trajectory $\Exp (\lambda, s)$ is denoted by
$$\tconj(\lambda)=\min \left\{t>0 \mid t \text{ is a conjugate time along } \Exp (\lambda, s), \ s \geq 0 \right\}. $$

The trajectory $\Exp (\lambda, s)$ loses its local optimality at the instant $t=\tconj(\lambda)$ (see~\cite{notes}).

The following lower bound on the first conjugate time is the main result of work~\cite{engel_conj}.
\begin{thrm}[\cite{engel_conj}]\label{th:tconjmax}
For any $\lambda \in C$
\begin{align}\label{tconjmax}
\tconj(\lambda) \geq \tmax(\lambda).
\end{align} 
\end{thrm}

\section{Decompositions in preimage and image of exponential mapping}
In this section we describe decomposition~(\ref{MtildeMi}) in the image, and decomposition~(\ref{NtildeDi}) in the preimage of the exponential mapping, which will be proved to be diffeomorphic via the exponential mapping in Th.~\ref{Expdif}.
\subsection{Decomposition in $M$} \label{decompM}

Let $\widehat{M} = M \backslash \{ q_0 \}$, then $M = \widehat{M} \sqcup \{ q_0 \}$. Further, we denote the subset containing the Maxwell strata $\MAX^1$ and $\MAX^2$:
$$ M' = \{ q \in \widehat{M} ~ | ~ xz = 0 \} $$
and its complement
$$ \widetilde{M} = \{ q \in \widehat{M} ~ | ~ xz \ne 0 \}, $$
then
\begin{align}\label{Mhattilde'}
&\widehat{M} = \widetilde{M} \sqcup M'. 
\end{align} 
Denote the connected components of the set $\widetilde{M}$:
\begin{align}
& M_1 = \{ q \in M ~ | ~ x < 0, z > 0 \}, \label{M1} \\ 
& M_2 = \{ q \in M ~ | ~ x < 0, z < 0 \}, \label{M2} \\ 
& M_3 = \{ q \in M ~ | ~ x > 0, z < 0 \}, \label{M3} \\ 
& M_4 = \{ q \in M ~ | ~ x > 0, z > 0 \}, \label{M4} 
\end{align}
so that
\begin{align}
&\widetilde{M} = \bigsqcup_{i=1}^{4}M_i. \label{MtildeMi}
\end{align}
This decomposition agrees with the action of reflections and dilations as described in the following statement.

\begin{prpstn} ~ \label{epsjdMi}

\begin{itemize}
\item[$(1)$] Reflections ${\varepsilon}^j \in G$ permute the domains $M_i$ according to Table~\ref{epsjMi}.
\item[$(2)$] Dilations $\delta_\mu, \ \mu > 0$, preserve the domains $M_i$.
\end{itemize}

\begin{table}[ht]

\begin{center}
\begin{tabular}{|c|c|c|c|}
\hline
$Id, {\varepsilon}^6$ & ${\varepsilon}^1, {\varepsilon}^7$ & ${\varepsilon}^2, {\varepsilon}^4$ & ${\varepsilon}^3, {\varepsilon}^5$ \\
\hline
$M_1$ & $M_2$ &$M_4$ & $M_3$ \\
\hline
$M_2$ & $M_1$ &$M_3$ & $M_4$ \\
\hline
$M_3$ & $M_4$ &$M_2$ & $M_1$ \\
\hline
$M_4$ & $M_3$ &$M_1$ & $M_2$ \\
\hline
\end{tabular}
\end{center}
\caption{Action of the reflections $\varepsilon^j$ on the domains $M_i$} \label{epsjMi}
\end{table}
\end{prpstn}

\begin{proof} Follows immediately from the definitions of the actions of reflections ${\varepsilon}^j \colon M \to M$, see~(\ref{epsiM})--(\ref{eps7M}),
and dilations $\delta_\mu \colon M \to M$, see~(\ref{dmuM}).
\end{proof}

\subsection{Decomposition in $N$}

Denote the subset in preimage of the exponential mapping that corresponds to all potentially optimal geodesics:
$$ \widehat{N} = \{ (\lambda, t) \in N ~ | ~ t \leq \tmax (\lambda) \}. $$

If $(\lambda, t) \in N \backslash \widehat{N}$, then the geodesic $\Exp(\lambda, s), \ s \in [0, t]$, is non-optimal. We decompose the set $\widehat{N}$ into subsets corresponding to the subsets of the set $\widehat{M}$ (Subsec.~\ref{decompM}),
the proof of this correspondence will be given in Subsec.~\ref{basicExp}. Let 
\begin{align*}
 &       N' = \{ (\lambda, t) \in N ~ | ~ t = \tmax (\lambda) \text{ or } \ c_{t/2} \sin \theta_{t/2} = 0 \}, \\
 &\widetilde{N} = \{ (\lambda, t) \in N ~ | ~ t < \tmax (\lambda), \ c_{t/2} \sin \theta_{t/2} \ne 0 \}, 
\end{align*}
then
\begin{align}\label{Nhattilde'} 
&\widehat{N} = \widetilde{N} \sqcup N'.  
\end{align}

\begin{figure}[htbp]
\centering
\includegraphics[width=0.32\linewidth]{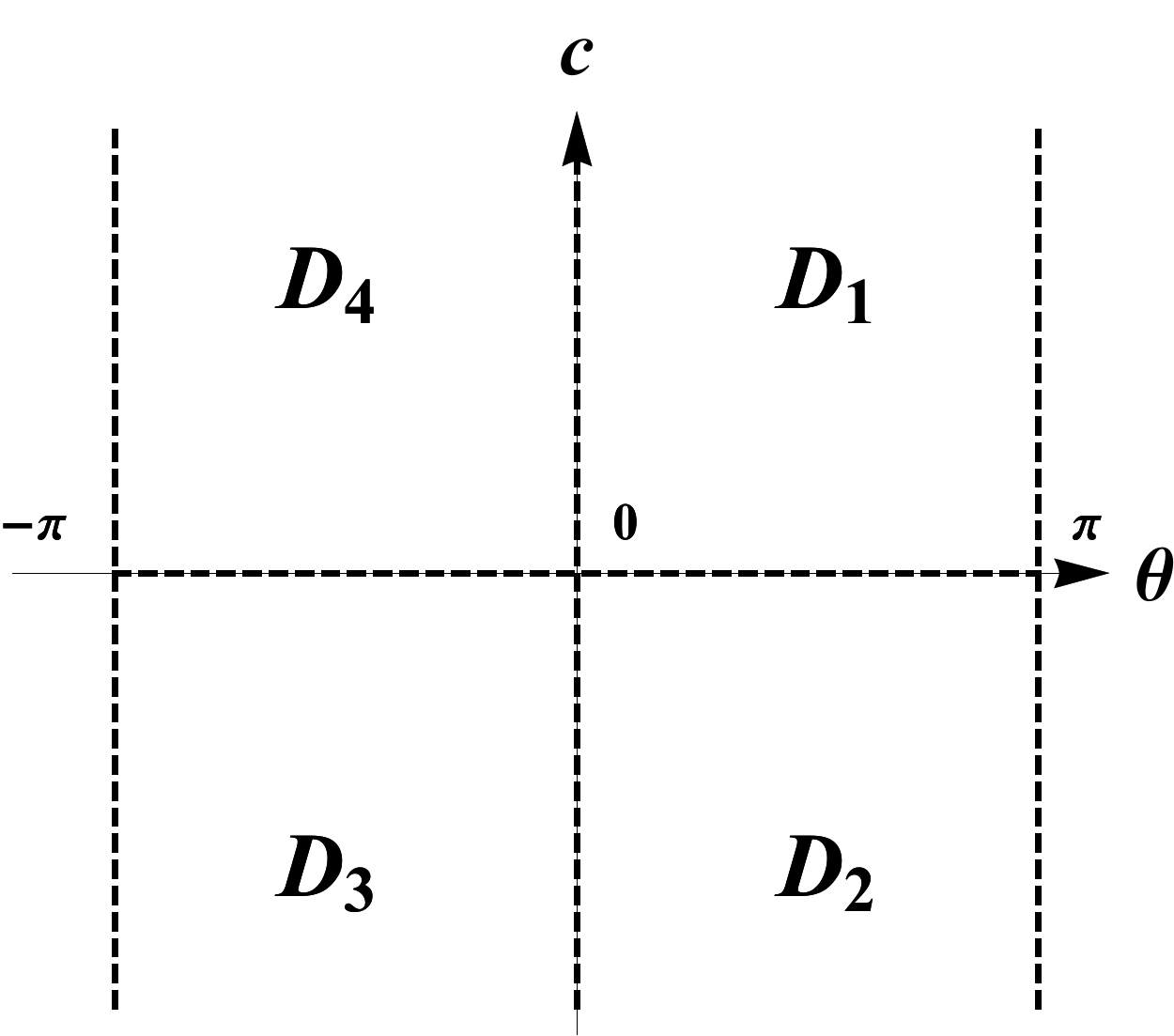}
\qquad
\includegraphics[width=0.3\linewidth]{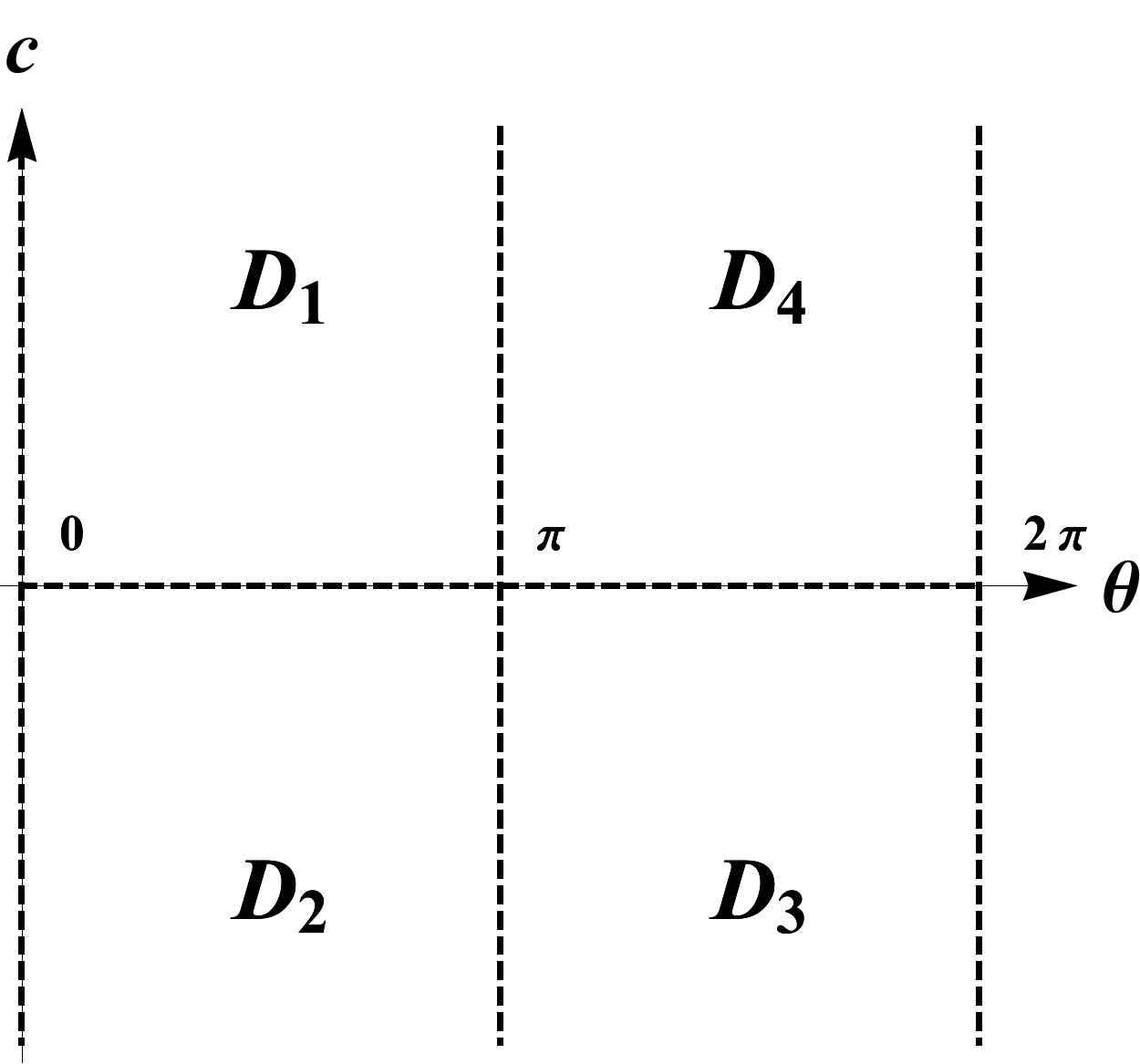}
\caption{The trace of domains $D_i$ in the set $\{t=0\}$ for $\alpha \geq 0$ and for $\alpha < 0$}\label{fig:Di}
\end{figure}

The following sets will play an important role in the description of the global structure of the exponential mapping:

\begin{align}
& D_1 = \{ (\lambda, t) \in N ~ | ~ t \in \big(0, \tmax(\lambda)\big), \ \sin\theta_{t/2} > 0, \ c_{t/2} > 0 \}, \label{D1}\\
& D_2 = \{ (\lambda, t) \in N ~ | ~ t \in \big(0, \tmax(\lambda)\big), \ \sin\theta_{t/2} > 0, \ c_{t/2} < 0 \}, \label{D2} \\
& D_3 = \{ (\lambda, t) \in N ~ | ~ t \in \big(0, \tmax(\lambda)\big), \ \sin\theta_{t/2} < 0, \ c_{t/2} < 0 \}, \label{D3} \\
& D_4 = \{ (\lambda, t) \in N ~ | ~ t \in \big(0, \tmax(\lambda)\big), \ \sin\theta_{t/2} < 0, \ c_{t/2} > 0 \}. \label{D4} 
\end{align}

We have the obvious decomposition
\begin{align} \label{NtildeDi}
 &\widetilde{N} = \bigsqcup_{i=1}^{4}D_i. 
\end{align}

The trace of domains $D_i$ in the set $\{(\lambda,t) \in N ~ | ~ t = 0\}$ is shown in Fig.~\ref{fig:Di}.

\begin{prpstn} ~\label{epsjdDi}

\begin{enumerate}
\item\label{pD1} Reflections ${\varepsilon}^j \in G$ permute the sets $D_i$ according to Table~\ref{epsjDi}. 
\item\label{pD2} Dilations $\delta_\mu , \mu > 0,$ preserve the sets $D_i$.
\end{enumerate}

\begin{table}[ht]

\begin{center}
\begin{tabular}{|c|c|c|c|}
\hline
$Id, {\varepsilon}^6$ & ${\varepsilon}^1, {\varepsilon}^7$ & ${\varepsilon}^2, {\varepsilon}^4$ & ${\varepsilon}^3, {\varepsilon}^5$ \\
\hline
$D_1$ & $D_2$ &$D_4$ & $D_3$ \\
\hline
$D_2$ & $D_1$ &$D_3$ & $D_4$ \\
\hline
$D_3$ & $D_4$ &$D_2$ & $D_1$ \\
\hline
$D_4$ & $D_3$ &$D_1$ & $D_2$ \\
\hline
\end{tabular}
\end{center}
\caption{Action of the reflections ${\varepsilon}^j$ on the domains $D_i$} \label{epsjDi}
\end{table}
\end{prpstn}

\begin{proof}
(\ref{pD1}) We prove only the equality ${\varepsilon}^1 (D_1) = D_2$, all the rest equalities given in Table~\ref{epsjDi} are proved similarly.
Let $(\lambda, t) = (\theta, c, \alpha, t) \in D_1$ and 
${\varepsilon}^1 (\lambda, t) = ({\lambda}^1 , t) = ({\theta}^1 , c^1, {\alpha}^1 , t),$ we show that $({\lambda}^1 , t) \in D_2$. 

Denote ${\lambda}_{t/2} = ({\theta}_{t/2} , c_{t/2} , \alpha) = e^{(t/2) \vec{H}_v} (\lambda)$ and ${\lambda}_{t/2}^1 = ({\theta}^1_{t/2} , c^1_{t/2}, {\alpha}^1) = e^{ (t/2) \vec{H}_v} ({\lambda}^1)$. Since ${\varepsilon}^1_* \vec{H}_v = - \vec{H}_v$, then ${\lambda}^1 = \varepsilon^1 \circ e^{t \vec{H}_v} (\lambda)$, thus
\begin{align*}
&{\lambda}^1_{t/2} = e^{(t/2) \vec{H}_v} \circ \varepsilon^1 \circ e^{t \vec{H}_v} (\lambda) = \varepsilon^1 \circ  e^{-(t/2) \vec{H}_v} \circ e^{t \vec{H}_v} (\lambda) = \varepsilon^1 \circ e^{(t/2) \vec{H}_v} (\lambda) = \varepsilon^1 ({\lambda}_{t/2}). 
\end{align*}
That is, $$({\theta}^1_{t/2}, c^1_{t/2}, {\alpha}^1) = ({\theta}_{t/2}, -c_{t/2}, \alpha).$$
The inclusion $(\lambda, t) \in D_1$ means that $$t \in \big(0, \tmax (\lambda)\big), \qquad \sin {\theta}_{t/2} > 0, \qquad c_{t/2} > 0,$$ thus $\sin {\theta}_{t/2}^1 > 0, \ c_{t/2}^1 < 0$. Moreover, since $\tmax ({\lambda}^1) = \tmax \circ {\varepsilon}^1 (\lambda) = \tmax (\lambda)$ by Propos.~\ref{tMaxinvar}, then $t \in \big(0, \tmax ({\lambda}^1)\big).$ Consequently, $({\lambda}^1 , t) \in D_2 .$

We proved that ${\varepsilon}^1 (D_1) \subset D_2.$ It follows similarly that ${\varepsilon}^1 (D_2) \subset D_1 $. Since ${\varepsilon}^1 \circ {\varepsilon}^1 = \mathrm{Id}$ on $N$, we have ${\varepsilon}^1 (D_1) = D_2$.

\medskip (\ref{pD2}) Let $(\lambda, t) = (\theta, c, \alpha, t) \in N, \ ({\lambda}_\mu , t_\mu) = \delta_\mu (\lambda, t) = (\theta, c/\mu, \alpha/\mu^2, \mu t). $ Since $\tmax ({\lambda}_\mu) = \mu \tmax (\lambda)$ by Propos.~\ref{tMaxinvar}, it is obvious that if $(\lambda, t) \in D_1$, then $({\lambda}_\mu, t_\mu) \in D_1$. Thus $\delta_\mu (D_1) \subset D_1.$ Since $d_{1/\mu} = (\delta_\mu)^{-1}$, then $\delta_\mu (D_1) = D_1$. It follows similarly that $\delta_\mu (D_i) = D_i$ for $i = 2, 3, 4.$
\end{proof}

\subsection{Basic properties of exponential mapping}\label{basicExp}
In this subsection we describe some simple properties on the action of the exponential mapping on the subsets of $N$ defined in the previous subsection.

First of all, $\Exp(\widehat{N}) \supset \widehat{M}$ since for any point $q_1 \in \widehat{M}$ there exists an optimal trajectory $q_s = \Exp(\lambda, s)$ such that $q_{t_1} = q_1$, thus $t_1 \leq \tcut (\lambda) \leq \tmax (\lambda)$, \ie, $\Exp(\lambda, t_1) = q_1$ with $(\lambda, t_1) \in \widehat{N}$. However, Maxwell points in $\widehat{M}$ have several preimages in $\widehat{N}$. Moreover, the mapping $\Exp|_{\widehat{N}}$ is degenerate
at points $(\lambda, t)$ where $t = \tmax(\lambda)$ is a conjugate time along the trajectory $\Exp(\lambda, s)$.

In the next two propositions we show that the action of $\Exp$ is compatible with decompositions~(\ref{Nhattilde'}), (\ref{Mhattilde'}), and (\ref{NtildeDi}), (\ref{MtildeMi}).

\begin{prpstn} \label{ExpN'}
There holds the inclusion 
\begin{align} \label{ExpN'xz=0}
& \Exp (N') \subset \{ q \in M ~ | ~ xz = 0 \} = M' \sqcup \{ q_0 \}.
\end{align}
\end{prpstn}

\begin{proof} The reflection ${\varepsilon}^4$ and dilations $\delta_\mu, \ \mu > 0,$ are symmetries of $\Exp$ and preserve the sets $N' , M'$ and $\{ q_0 \}$. Since ${\varepsilon}^4\colon \alpha \mapsto -\alpha$ and $\delta_\mu \colon \alpha \mapsto \alpha / \mu^2$, we can assume in the proof of inclusion~(\ref{ExpN'xz=0}) that $\alpha \in \{ 0, 1 \}$.

Let $(\lambda, t) \in N'$ and $q_t = (x_t, y_t, z_t, v_t) = \Exp(\lambda, t)$, we show that $x_t z_t = 0$.

Suppose first that $\alpha=1$, then $\lambda \in \bigcup_{i=1}^{5}C_i$.

Let $\lambda \in C_1$, then we use parameterization of extremals~(\ref{ExpC1}). Since $(\lambda, t) \in N'$, then $c_{t/2} \sin \theta_{t/2} = 0$ or $t = \tmax (\lambda)$.
If $c_{t/2} = 2 k \cn \tau = 0$, then $\cn \tau = 0$, thus $z_t = 0$ in view of~(7.3)~\cite{engel}. If $\sin {\theta}_{t/2} = 2 k \sn \tau \dn \tau = 0,$ then $x_t = 0$ in view of~(7.2)~\cite{engel}. Finally, if $t = \tmax (\lambda)$, then $p = p_z^1(k)$ or $p = 2 K(k)$ by~(\ref{tmaxC1}), thus $z_t = 0$ or $x_t = 0$ by~(7.2) and (7.3)~\cite{engel}. 

The case $\lambda \in C_2 \cup C_3$ is considered similarly to the case $\lambda \in C_1$.

If $\lambda \in C_4 \cup C_5$, then $x_t = 0$ by~(\ref{ExpC4}),~(\ref{ExpC5}).

Now suppose that $\alpha = 0$, thus $\lambda \in C_6 \cup C_7$.

Let $\lambda \in C_6$, then we use parameterization of extremals~(\ref{ExpC6}). The case $c_{t/2} = 0$ is impossible. If $\sin {\theta}_{t/2} = 0$, then $t = 2\pi /|c|$, thus $x_t = 0$. If $t = \tmax (\lambda) = 2\pi /|c|~(\ref{tmaxC6})$, then $x_t = 0$ as well.

Finally, if $\lambda \in C_7,$ then $z_t = 0$ by~(\ref{ExpC6}).
\end{proof}

\begin{prpstn} \label{ExpDiMi}
For any $i = 1, \ldots, 4$, we have $Exp(D_i) \subset M_i$.
\end{prpstn}
\begin{proof} By virtue of the reflections ${\varepsilon}^i$ (Propos.~\ref{epsjdMi},~\ref{epsjdDi}), the proof of this proposition reduces to the case $i = 1.$ So let $(\lambda, t) \in D_1$, we prove that $q_t = \Exp(\lambda, t) \in M_1$.

The reflection ${\varepsilon}^6$ and the dilations $\delta_\mu , \mu > 0$, preserve the domains $D_1$ and $M_1$, and act on the parameter $\alpha$ as $\ds {\varepsilon}^6 : \alpha \mapsto -\alpha, \ \delta_\mu : \alpha \mapsto \frac{\alpha}{\mu^2}$, thus we can assume in this proof that $\alpha \in \{ 0, 1 \}$.

Since $(\lambda, t) \in D_1$, then $\sin {\theta}_{t/2} > 0, \ c_{t/2} > 0, \ t \in \big(0, \tmax (\lambda)\big).$

Let $\alpha = 1, $ then $\lambda \in C_1 \cup C_2 \cup C_3.$

Let $\lambda \in C_1.$ Then $\sin {\theta}_{t/2} = 2 k \sn \tau \dn \tau >0, \ c_{t/2} = 2 k \cn \tau > 0.$ Since in this case $t \in \big(0, \tmax (\lambda)\big)$ and $\tmax (\lambda) = \min (2 p_z^1, 4 K),$ then $f_z (p, k) > 0$ and $\sn p > 0.$ Then formulas~(7.2),~(7.3)~\cite{engel} imply that $x_t < 0, z_t > 0$, \ie, $q_t \in M_1$.

The cases $\lambda \in C_2$ and $\lambda \in C_3$ are considered similarly to the case $\lambda \in C_1$.

Now let $\alpha = 0,$ then $\lambda \in C_6$. Then $\ds x_t = - \frac{2}{c} \sin {\theta}_{t/2} \sin \frac{c t}{2} < 0$ and $\ds z_t = \frac{c t - \sin (c t)}{2 c^2} > 0$, thus $q_t \in M_1$.

We proved that $\Exp(D_1) \subset M_1$.
\end{proof}

Our goal is to prove that the mappings $\Exp \colon D_i \to M_i, \ i = 1, \dots, 4$, are diffeomorphisms, see Th.~\ref{Expdiffeo}. This is done in Section~\ref{sec:Dif} via the following Hadamard global diffeomorphism theorem.

\begin{thrm}[\cite{Hadamard}] \label{Hadamard}
Let $F \colon X \to Y$ be a smooth mapping between smooth manifolds of equal dimension. Let the following conditions hold:
\begin{enumerate}
\item [$(1)$] $X$ is connected,
\item [$(2)$] $Y$ is connected and simply connected,
\item [$(3)$] $F$ is nondegenerate,
\item [$(4)$] $F$ is proper \big(\ie, $F^{-1} (K) \subset X$ is compact for a compact $K \subset Y$\big).
\end{enumerate}
Then $F$ is a diffeomorphism.
\end{thrm}

\subsection{Topological properties of decompositions in $M$ and $N$}
We prove that hypotheses (1), (2) of Th.~\ref{Hadamard} are verified for the mappings $\Exp \colon D_i \to M_i$.

\begin{dfntn} Suppose $X$ is a topological space and $f_1, f_2 \colon X \to \xR$. Then $f_1 \sim f_2$ on a sequence $\{\lambda_n\} \subset X$ if $\ds \lim_{n\rightarrow\infty} \frac{f_1 (\lambda_n)}{f_2(\lambda_n)} = 1$. 
\end{dfntn}

\begin{prpstn} ~\label{DiMitopo}

\begin{enumerate}
\item \label{pDM1} The sets $D_i \subset N, \ i = 1,\dots,4,$ are open and connected.
\item \label{pDM2} The sets $M_i \subset M, \ i = 1,\dots,4,$ are open, connected and simply connected.
\end{enumerate}
\end{prpstn}

In the proof of item (\ref{pDM1}) of this proposition we need the following statement.

\begin{prpstn} \label{tMaxcont}
The function $\tmax \colon C \to (0, +\infty ]$ is continuous on the set $C \backslash C_4$, and is smooth on the set $C^0_1 \cup C_2,$ where $C^0_1 = \{ \lambda \in C_1 ~ | ~ k \ne k_0 \}$.
\end{prpstn}

\begin{rmrk} 
We assume in $(0, +\infty ]$ the natural basis of topology:
$$ (a, b), \qquad (a, +\infty], \quad 0 < a < b < +\infty.$$
\end{rmrk}

\begin{proof}
Let $ {\lambda}_n \to \bar{\lambda}$ as $n \to \infty$, where ${\lambda}_n , \bar{\lambda} \in C \backslash C_4 = (\cup^{3}_{i=1} C_i) \cup (\cup^{7}_{i=5} C_i)$. We denote $t_n = \tmax ({\lambda}_n)$ and $\bar{t} = \tmax (\bar{\lambda})$, then prove that $t_n \to \bar{t}$ as $n \to +\infty$. 

\medskip 1. Let ${\lambda}_n \in C_1,$ then $\bar{\lambda} \in \cl(C_1) \backslash C_4 = C_1 \cup C_3 \cup C_5 \cup C_7.$ 

\medskip 1.1. Let $\bar{\lambda} \in C_1.$ The function  $ \ds \tmax |_{_{C_1}} = \frac{\min\big(2 p_z^1 (k), 4 K(k)\big)}{\sigma}$
				is continuous since for $k \in (0, 1)$ the function $\min\big(p_z^1 (k), 2 K(k)\big)$ is continuous (see Cor.~3.1~\cite{el_exp}), thus $t_n \to \bar{t}$. 

\medskip 1.2. Let $\bar{\lambda} \in C_3 \cup C_5$. Then $k_n = k({\lambda}_n) \to 1, \ K(k_n) \to +\infty, \ p_z^1 (k_n) \to +\infty, \ \sigma(\lambda_n) \to \bar{\sigma} > 0$. Thus $t_n \to +\infty = \bar{t}$. 

\medskip 1.3. Let $\bar{\lambda} \in C_7$. Then ${\alpha}_n \to 0$. Since $\min\big(2 p_z^1 (k), 4 K(k)\big) > 2 K (k) > \pi$, then $t_n \to +\infty = \bar{t}$. 

\medskip 2. Let ${\lambda}_n \in C_2$, then $\bar{\lambda} \in \cl (C_2) \backslash C_4 = C_2 \cup C_3 \cup C_5 \cup C_6 \cup C_7.$ 

\medskip 2.1. Let $\bar{\lambda} \in C_2$. The function $\ds \tmax |_{_{C_2}} = \frac{2 K(k)k}{\sigma} $ is continuous, thus $t_n \to \bar{t}$. 

\medskip 2.2. Let $\bar{\lambda} \in C_3 \cup C_5.$ This case is similar to case 1.2. 

\medskip 2.3. Let $\bar{\lambda} \in C_6$. Then 
			\begin{align*}
				&{\alpha}_n = \alpha({\lambda}_n) \to \bar{\alpha} = 0, &&c_n = c({\lambda}_n) \to \bar{c} \ne 0, \\
				&E_n = E ({\lambda}_n) \to \frac{\bar{c}^2}{2} = \bar{E} \ne 0, &&k_n = k({\lambda}_n) = \sqrt{ \frac{2 |{\alpha}_n|}{ E_n + |{\alpha}_n|} } \sim 
									\frac{2 \sqrt{ |{\alpha}_n|} }{ |\bar{c}| } \to 0, \\
				&t_n = \frac{2 K (k_n) k_n}{ \sqrt{|{\alpha}_n|} } \sim 2 K (0) \cdot \frac{2}{|\bar{c}|} = \frac{2 \pi}{|\bar{c}|} = \bar{t},&&
			\end{align*}
 \ie, $t_n \to \bar{t}$.
 
\medskip 2.4. Let $\bar{\lambda} \in C_7$. Then ${\alpha}_n \to 0, c_n \to 0.$ Thus $E_n \to 0,$ so $\ds \frac{k_n}{ \sqrt{|{\alpha}_n|} } = \sqrt{ \frac{2}{E_n + |{\alpha}_n|} } \to +\infty$. 
					Consequently, $\ds t_n = \frac{2 K (k_n) k_n }{ \sqrt{ |{\alpha}_n| } } \to +\infty = \bar{t}$. 

\medskip 3. Let ${\lambda}_n \in C_3$, then $\bar{\lambda} \in \cl(C_3) \backslash C_4 = C_3 \cup C_5 \cup C_7,$ and this case is similar to cases 1.2 and 1.3. 

\medskip 4. Let ${\lambda}_n \in C_5$, then $\bar{\lambda} \in \cl(C_5) = C_5 \cup C_7$, and $t_n = +\infty = \bar{t}$. 

\medskip 5. Let ${\lambda}_n \in C_6$, then $\bar{\lambda} \in \cl(C_6) = C_6 \cup C_7$. 

\medskip 5.1. Let $\bar{\lambda} \in C_6$. Since the function $\ds \tmax |_{_{C_6}} = \frac{2 \pi}{|c|}$
					is continuous, then $t_n \to \bar{t}$. 

\medskip 5.2. Let $\bar{\lambda} \in C_7$. Then $c_n \to 0$, thus $\ds t_n = \frac{2 \pi}{ |c_n| } \to +\infty = \bar{t}$.
				
\bigskip The function $\tmax (\lambda)$ is smooth on $C^0_1$ since for $\lambda \in C_1$ we have by virtue of~(\ref{tmaxC1}):
					\begin{align*}
						&k < k_0 ~ \Rightarrow ~ \tmax = \frac{ 2 p_z^1 (k) }{ \sqrt{ |\alpha| } } \in C^{\infty}, \\
						&k > k_0 ~ \Rightarrow ~ \tmax = \frac{4 K }{ \sqrt{ |\alpha| } } \in C^{\infty}. 
					\end{align*}
					Similarly, $\tmax (\lambda)$ is smooth on $C_2$ by virtue of~(\ref{tmaxC2}). The proof of Propos.~\ref{tMaxcont} is complete.
\end{proof}

\begin{rmrk} 
The function $\tmax$ is discontinuous on $C_4$.
\end{rmrk}

Indeed, let ${\lambda}_n \in C_1$ be such that $k ({\lambda}_n) \to 0$ and $\alpha ({\lambda}_n) \to \bar{\alpha} \ne 0. $ Then ${\lambda}_n \to \bar{\lambda} \in C_4$ but 
$$ \tmax ({\lambda}_n) \to \frac{ 2 p_z^1 (0) }{ \sqrt{ |\bar{\alpha}| } } < \tmax (\bar{\lambda}) = +\infty. $$
Here $p = p_z^1 (0)$ is the minimal positive root of the equation $f_z (p, 0) = \sin p - p \cos p = 0$, thus $p \in (\pi, 3 \pi/2)$. 

Now we prove Propos.~\ref{DiMitopo}.
\begin{proof} (\ref{pDM1}) Reflections ${\varepsilon}^i \colon N \to N$ are diffeomorphisms and permute the sets $D_i$, thus it is sufficient to prove that the set $D_1 = \{ (\lambda, t) \in N ~ | ~ \sin {\theta}_{t/2} > 0, \ c_{t/2} > 0, \ t < \tmax (\lambda)\}$ is open and connected.

Consider the vector field $P = \frac{t}{2} (c \frac{\partial}{\partial \theta} - \alpha \sin \theta \frac{\partial}{\partial c}) \in \operatorname{Vec}(N)$. Denote the flow of this vector field for time 1 as
$e^P \in \Diff(N)$. We have
$$ e^{P} (\theta, c, \alpha, t) = e^{P} (\lambda, t) = \big(e^{ \frac{t}{2} \vec{H}_v } (\lambda), t \big) = ({\theta}_ {t/2}, c_{t/2}, \alpha, t), $$
thus $e^{P} (D_1) = \widetilde{D}_1$, where

$$ \widetilde{D}_1 = \{ (\lambda , t) \in N ~ | ~ \sin \theta > 0, \ c > 0, \ t < \tmax (\lambda, t) \} .$$

By Propos.~\ref{tMaxcont}, the function $\big(t - \tmax (\lambda)\big) \colon N ~ \to ~ (0, +\infty]$ is continuous on the set $N \backslash N_4 \supset \widetilde{D}_1,$ thus the set $\widetilde{D}_1$ is open. Moreover, the domain $\widetilde{D}_1$ is a subgraph of the function $\tmax (\lambda)$ on a connected domain $\{ (\theta, c, \alpha) \in C ~ | ~ \theta \in (0, \pi), \ c>0, \ \alpha \in \mathbb{R} \}$, thus $\widetilde{D}_1$ is connected. 

We proved that $\widetilde{D}_1$ is open and connected, thus $D_1 =  e^{ - P} (\widetilde{D}_1)$ is open and connected as well.

\medskip (\ref{pDM2}) It is obvious from definitions~(\ref{M1})--(\ref{M4}) that the sets $M_i, \ i = 1, \dots 4$, are open, connected and simply connected.
\end{proof}

\section{Diffeomorphic properties of exponential mapping}\label{sec:Dif}
In this section we prove that restriction of the exponential mapping to the subdomains $D_i, M_i$ is a diffeomorphism.

\begin{lmm}\label{lm1}
If $\, \Exp : D_1 \rightarrow M_1$ is proper, then $\Exp : D_i \rightarrow M_i$ is proper for $i=2,3,4$.
\end{lmm}
\begin{proof}
Follows immediately from Propositions~\ref{epsjdMi},~\ref{epsjdDi}.
\end{proof}

\begin{lmm}\label{lm2}
The mapping $\Exp : D_1 \rightarrow M_1$ is proper iff there exists no sequence  
$\{\nu_n \} \subset \newD = (D_1 \cap N_1) \cup (D_1 \cap N_{2}) \cup (D_1 \cap N_{3}) \cup (D_1 \cap N_{6})$, such that $\nu_n \rightarrow \bar{\nu} \in \cl (D_1) \backslash D_1$ and $\Exp(\nu_n) \rightarrow \bar{q} \in M_1$. 
\end{lmm}
\begin{proof}
It follows from the definition of a proper mapping that the mapping $\Exp : D_1 \rightarrow M_1$ is proper iff there exists no sequence $\{\nu_n \} \subset D_1$, such that $\nu_n \rightarrow \bar{\nu} \in \cl (D_1) \backslash D_1$ and $\Exp(\nu_n) \rightarrow \bar{q} \in M_1$.

Moreover, the definition of $D_1$~(\ref{D1}) gives the decomposition 
$$D_1 = (D_1 \cap N_1) \cup (D_1 \cap N_{2}) \cup (D_1 \cap N_{3}) \cup (D_1 \cap N_{6}).$$
\end{proof}

Let us introduce the following sets for arbitrary $\varepsilon \in (0,1)$: 
$$\ds S_\varepsilon := \{\nu \in N ~ | ~ \theta_{t/2} \in [\varepsilon, \pi-\varepsilon], \ c_{t/2} \in[\varepsilon, 1/\varepsilon], \ |\alpha| \leq 1/\varepsilon, \ t \in [\varepsilon, 1/\varepsilon], \ \tmax(\lambda)-t \geq \varepsilon\}.$$ 

\begin{lmm}
The set $S_\varepsilon$ is compact for any $\varepsilon > 0$.  
\end{lmm}
\begin{proof}
Let $\{\nu_n\} \subset S_\varepsilon$ be an arbitrary sequence. To prove the lemma, we need to find a subsequence $\nu_{n_m}$ which tends to $\bar{\nu} \in S_\varepsilon$ as $m \rightarrow \infty$. 

Since $\alpha, t$ are bounded on $S_\varepsilon$, we obtain for a subsequence that $\alpha \rightarrow \bar{\alpha}, t \rightarrow \bar{t}$ as $m \to \infty$. 

Since $\theta_{t/2}, c_{t/2}$ are bounded on $S_\varepsilon$, we obtain for a subsequence $(\theta_{t/2}, c_{t/2}) \rightarrow (a,b)$. Moreover, we have 
$(\theta, c) = \Phi_{-t/2} (\theta_{t/2},c_{t/2}) \rightarrow \Phi_{-\bar{t}/2} (a,b) =: (\bar{\theta}, \bar{c}),$ 
where $\Phi$ is the flow of pendulum~(\ref{pend}). Since $\bar{\nu}= (\bar{\theta}, \bar{c}, \bar{\alpha}, \bar{t}) \in S_\varepsilon$ by continuity of the functions which define $S_\varepsilon$, we see that $S_\varepsilon$ is compact. 
\end{proof}

\begin{lmm}
If $K \subset D_1$ is compact, then there exists $\varepsilon>0$ such that $K \subset S_\varepsilon$.
\end{lmm}

\begin{proof}
Since the functions $\theta_{t/2}, \, c_{t/2}, \, \alpha, \, t, \, (\tmax - t) \, $ are continuous on $N$, these functions attain maximum and minimum on $K$.
\end{proof}

\begin{lmm}
Let $\{\nu_n\} \subset D_1$. Then  $\nu_n \rightarrow  \bar{\nu} \in \cl(D_1) \backslash D_1$ iff one of the following conditions holds for $\{\nu_n\}$:
\begin{enumerate}
\item \label{cond1} $\theta_{t/2} \rightarrow 0$,
\item $\theta_{t/2} \rightarrow \pi$,
\item $c_{t/2} \rightarrow 0$,
\item $c_{t/2} \rightarrow +\infty$,
\item $t \rightarrow 0$,
\item $\tmax (\lambda) - t \rightarrow 0$,
\item \label{cond7} $|\alpha| \rightarrow \infty$.
\end{enumerate}
\end{lmm}
\begin{proof}Necessity. Assume the converse. Suppose for any sequence $\{\nu_n\} \subset D_1, \ \nu_n \rightarrow \bar{\nu} \in \cl(D_1) \backslash D_1$, that conditions~(\ref{cond1})--(\ref{cond7}) do not hold. This means that there exists $\varepsilon > 0$ such that conditions 
$$\theta_{t/2} \geq \varepsilon, \ \theta_{t/2} \leq \pi - \varepsilon, \quad c_{t/2} \geq \varepsilon, \ c_{t/2} \leq 1/\varepsilon, \quad t \geq \varepsilon, \ \tmax (\lambda) - t \geq \varepsilon, \quad  |\alpha| \leq 1/\varepsilon$$ 
hold for a subsequence. It follows that $\{\nu_n\} \subset S_\varepsilon$, which is a compact subset of $D_1$. So $\bar{\nu} \in S_\varepsilon \subset D_1$. This contradiction proves the necessity.

Sufficiency. Assume the converse. Let for any sequence $\nu_n \subset D_1$ we have $\nu_n \rightarrow \bar{\nu} \in D_1$. Then there exists a compact set $K \supset \{\nu_n\}, \ \bar{\nu} \in K$. This means that there exists $\varepsilon > 0$ such that $K \subset S_\varepsilon$. This contradiction proves the lemma.
\end{proof}

\begin{dfntn} Suppose $X$ is a topological space and $f_1, f_2 \colon X \to \xR$. Then $f_1 \newtilde f_2$ on a sequence $\{\nu_n\} \subset X$ if $\ds \lim_{n\rightarrow\infty} \frac{f_1 (\nu_n)}{f_2(\nu_n)} \in \xR \backslash \{0\}$. 
\end{dfntn}

In the next lemmas we use the parametrization of exponential mapping for the case $\lambda \in C_6$, see~(\ref{ExpC6}).

\begin{lmm} \label{N61}
If $\{\nu_n\} \subset D_1 \cap N_{6}, \ \nu_n \rightarrow \bar{\nu} \in \cl(D_1) \backslash D_1,$ and $\Exp (\nu_n) \rightarrow \bar{q}  \in M_1$, then $c \rightarrow 0$ on the sequence~$\{\nu_n\}$.
\end{lmm}
\begin{proof}
Notice that for $\bar{\nu}=\{\bar{x},\bar{y},\bar{z},\bar{v}\}\in M_1$ we have $\bar{x} \neq 0$.
Consider all possible cases $\bar{\nu} \in \cl{D_1} \backslash D_1$:

\medskip 1. $\ds \theta_{t/2} \rightarrow 0 ~ \Rightarrow ~ \frac{c t}{2} + \theta \rightarrow 0 \Rightarrow c\rightarrow 0 \text{ or } \ds x = -\frac{2 \sin (\frac{c t}{2} +\theta) \sin \frac{c t}{2}}{c} \rightarrow \bar{x} = 0$.

\medskip 2. $\ds \theta_{t/2} \rightarrow \pi ~ \Rightarrow ~ \frac{c t}{2} + \theta \rightarrow 0 \Rightarrow c\rightarrow 0 \text{ or } x \rightarrow \bar{x} = 0$.

\medskip 3. $\ds c_{t/2} \rightarrow 0 ~ \Rightarrow ~ c \rightarrow 0$.

\medskip 4. $\ds c_{t/2} \rightarrow \infty ~ \Rightarrow ~ c \rightarrow \infty \Rightarrow x \rightarrow \bar{x} = 0$.

\medskip 5. $\ds t\rightarrow 0 ~ \Rightarrow ~ c \rightarrow 0$, otherwise $x \rightarrow \bar{x} = 0$.

\medskip 6. $\ds t \rightarrow \frac{2 \pi}{|c|}$. This means that $c\rightarrow0$ or $x \rightarrow \bar{x} = 0$. 
\end{proof}

\begin{lmm} \label{N62}
Suppose $\nu_n \in D_1 \cap N_{6}$. If $c \rightarrow 0$, then $x \rightarrow 0$ or one of the functions $x,y,z$ tends to $\infty$ on the sequence $\{\nu_n\}$.
\end{lmm}
\begin{proof}
Consider two possible cases:

\medskip 1. If $c t \rightarrow 0$, then 
\begin{align*}
z &\newtilde \frac{(c t)^3}{c^2} = c t^3, \\
x^2+y^2 &= \frac{2 - 2\big(\cos(c t + \theta) \cos \theta + \sin(c t + \theta)\sin \theta\big)}{c} = \frac{2 \big(1- \cos (c t)\big)}{c} \newtilde c t^2. 
\end{align*}
It follows that $t \rightarrow \infty$, otherwise $x^2 \rightarrow 0$. Then we have $z \newtilde (x^2+y^2) t$. This means that $z \newtilde t \rightarrow \infty$ or $x^2 \rightarrow 0$.

\medskip 2. If $c t \rightarrow \bar{c} \neq 0 $ or $c t \rightarrow \infty$, then $c t - \sin (c t) > M > 0$, thus $z \rightarrow \infty$.
\end{proof}

In the next lemmas we use the following parametrization of exponential mapping for the case $\lambda \in C_3$ (see~(\ref{ExpC3})):
\begin{align*}
x&=-\frac{8 \sgn c\ \sigma \sh p \sh \tau}{\alpha \big(\ch(2 p) +\ch(2 \tau)\big)}, \\  
y&=\frac{2 \sigma}{\alpha} \bigg(\frac{2 \sh (2 p)}{\ch(2 p)+\ch(2 \tau)} - p \bigg), \\ 
z&=\frac{8 \sgn c\ \ch \tau (p \ch p -\sh p) } {|\alpha| (\ch(2 p)+\ch(2 \tau))}, \\
v&=- \frac{1}{3 \alpha \sigma \ch(p-\tau) \ch^2(p+\tau)} \Bigg(6 \big(\ch \tau - 3 \ch (2 p + \tau)\big) \sh p  \\
 &+  2 p \bigg(6 \ch (3 p + \tau) + p \ch (p + \tau) \Big(p \big(\ch (2 p) + \ch (2 \tau)\big) - 6 \sh (2 p)\Big)\bigg)\Bigg).             
\end{align*}

\begin{lmm} \label{N31}
If $\{\nu_n\} \subset D_1 \cap N_{3}$, $\nu_n \rightarrow \bar{\nu} \in \cl(D_1) \backslash D_1$, and $\Exp (\nu_n) \rightarrow \bar{q} \in M_1$, then $\sigma \rightarrow 0$ or $p \rightarrow \infty$ and $t \rightarrow \infty$ with $\sigma \rightarrow \bar{\sigma} \neq 0$. 
\end{lmm}
\begin{proof}
Notice that for $\bar{\nu}=\{\bar{x},\bar{y},\bar{z},\bar{v}\}\in M_1$ we have $\bar{x} \neq 0$. Consider all possible cases $\nu \rightarrow \cl(D_1) \backslash D_1$:

\medskip 1. $\ds \theta_{t/2} \rightarrow 0 ~ \Rightarrow ~ 
           \left\{  
           \begin{array}{rcl}  
           \tangh \tau \rightarrow 0\\  
           \ch \tau \rightarrow 1
           \end{array}   
           \right. ~\Rightarrow \tau \rightarrow 0 ~ \Rightarrow ~ \sigma \rightarrow 0$, or $x\rightarrow \bar{x} = 0$, or $z \rightarrow \infty$ with $p \rightarrow \infty$.

\medskip 2. $ \ds \theta_{t/2} \rightarrow \pi \Rightarrow 
           \left\{  
           \begin{array}{rcl}  
           \tangh \tau \rightarrow 1\\  
           \ds \frac{1}{\ch \tau} \rightarrow 0  
           \end{array}   
           \right. \Rightarrow  \tau \rightarrow \infty ~ \Rightarrow 
					 \left\{  
           \begin{array}{rcl}  
           p \rightarrow \infty \\  
           \tau \rightarrow \infty   
           \end{array}\right.   \text{or  } \sigma \rightarrow 0.$       

\medskip 3. $\ds c_{t/2} \rightarrow 0 ~ \Rightarrow ~ \frac{\sigma}{\ch \tau} \rightarrow 0 ~ \Rightarrow ~ 
           \left[  
           \begin{array}{rcl}  
           \sigma \rightarrow 0 \\  
           \ch \tau \rightarrow \infty   
           \end{array}\right. ~ \Rightarrow ~
           \left\{  
           \begin{array}{rcl}  
           p \rightarrow \infty \\  
           \tau \rightarrow \infty   
           \end{array}\right.   \text{ or  } \sigma \rightarrow 0.$

\medskip 4. $\ds c_{t/2} \rightarrow \infty ~ \Rightarrow ~ \frac{\sigma}{\ch \tau} \rightarrow \infty ~ \Rightarrow ~ \sigma \rightarrow \infty ~ \Rightarrow ~ x\rightarrow \bar{x} = 0.$

\medskip 5. $\ds t \rightarrow 0 ~ \Rightarrow ~ \frac{p}{\sigma} \rightarrow 0 ~ \Rightarrow ~ 
           \left[  
           \begin{array}{rcl}  
           p \rightarrow 0 \\  
           \sigma \rightarrow \infty   
           \end{array}\right. ~ \Rightarrow ~
           \left[  
           \begin{array}{rcl}  
           \sigma \rightarrow 0 \\  
           x \rightarrow 0   
           \end{array}\right. ~ \Rightarrow ~ \sigma \rightarrow 0.$

\medskip 6. $\ds t \rightarrow \infty ~ \Rightarrow ~ \frac{p}{\sigma} \rightarrow \infty ~ \Rightarrow ~
        \left[  
         \begin{array}{rcl}  
         p \rightarrow \infty \\  
          \sigma \rightarrow 0 
          \end{array}\right. ~ \Rightarrow ~
          \left\{  
          \begin{array}{rcl}  
          p \rightarrow \infty \\  
          \tau \rightarrow \infty   
          \end{array}\right.   \text{ or } \sigma \rightarrow 0.$

\medskip 7. $|\alpha| \rightarrow \infty ~ \Rightarrow ~ \left[  
           \begin{array}{rcl}  
           p \rightarrow \infty, \\  
           \tau \rightarrow \infty,   
           \end{array}\right.$ or $x \rightarrow \bar{x} = 0$        
\end{proof}

\begin{lmm} \label{N32}
Suppose $\nu_n \in D_1 \cap N_{3}$. If $p \rightarrow \infty, \tau \rightarrow \infty, \sigma \rightarrow \bar{\sigma} \neq 0$, then $y \rightarrow \infty$.
\end{lmm}
\begin{proof}
Since $\ds \frac{2 \sh (2 p)}{\ch (2p) + \ch(2\tau)} < \infty$ for $p \rightarrow \infty$, then $y \rightarrow \infty$.
\end{proof}

\begin{lmm} \label{N33}
Suppose $\nu_n \in D_1 \cap N_{3}$. If $\sigma \rightarrow 0$, then one of the functions $x,y,z,v$ tends to $\infty$, otherwise $x$ or $z$ tends to $0$.
\end{lmm}
\begin{proof}
Assume the converse. We have $\ds z \newtilde \frac{\ch \tau (p \ch p - \sh p)}{\sigma^2 \big(\ch (2 p) + \ch(2 \tau)\big)} ~ \Rightarrow ~
\left[  
         \begin{array}{rcl}  
         p \rightarrow \infty, \\  
         p \rightarrow 0, \\
         \tau \rightarrow \infty.
          \end{array}\right.$ 
So the proof is in these three cases as follows:

\medskip 1. $p \rightarrow \infty$. Then we get $\ds y \newtilde \frac{1}{\sigma} p \rightarrow \infty$.

\medskip 2. $p \rightarrow 0$. Consider three subcases:

\medskip 2.1. $\tau \rightarrow 0$. Here we have $\ds x \newtilde \frac{\tau p}{\sigma} \Rightarrow \sigma \newtilde \tau p$ \, and \, $\ds y \newtilde \frac{p}{\sigma} \newtilde \frac{p}{\tau p} = \frac{1}{\tau} \rightarrow \infty$.

\medskip 2.2. $\tau \rightarrow \infty$. We obtain 
\begin{align*}
&x \newtilde \frac{p}{\sigma e^\tau} \Rightarrow p \newtilde \sigma e^\tau, \\
&z \newtilde \frac{p^3}{\sigma^2 e^\tau} \newtilde \frac{p^2}{\sigma} \Rightarrow \sigma \newtilde p^2, \\
&y \newtilde \frac{1}{\sigma} \Big( \frac{p}{\ch^2 \tau} - p\Big) \newtilde \frac{p}{\sigma} \newtilde \frac{1}{p} \rightarrow \infty.
\end{align*}

\medskip 2.3. $\tau \rightarrow \bar{\tau} < \infty, \ \bar{\tau} \neq 0$. It follows that $\ds x \newtilde \frac{p}{\sigma} \Rightarrow p \newtilde \sigma$, then $\ds z \newtilde \frac{p^3}{\sigma^2} \newtilde \sigma \rightarrow 0$.

\medskip 3. $\tau \rightarrow \infty, \ p \rightarrow \bar{p} <\infty, \ \bar{p} \neq 0$. We get
$\ds \frac{\sh(2 p)}{\ch(2 p) + \ch(2 \tau)} \rightarrow 0 \Rightarrow y \newtilde \frac{\bar{p}}{\sigma} \rightarrow \infty$.
\end{proof} 

Below in the case $\nu_n \in (D_1 \cap N_{1}) \cup (D_1 \cap N_{2})$ we use the following notation:
\begin{align}
&s_{i} = \sin u_i, \ c_{i} = \cos u_i, \ d_{i} = \sqrt{1-k^2 s_{i}^2}, \quad i=1,2. \\
&E_1 = E(u_1, k), \ F_1 = F(u_1, k), \ \Delta = 1 - k^2 s_{1}^2 s_{2}^2.
\end{align}

\begin{lmm} \label{lemm1}
Suppose $\nu_n \in (D_1 \cap N_{1}) \cup (D_1 \cap N_{2})$. If $\Delta \rightarrow 0$, then $\ds \frac{d_{1} d_{2}}{\Delta}$ and $\ds \frac{c_{1} d_{1}}{\Delta}$ are bounded from above.
\end{lmm}
\begin{proof}
Let $k^2 = 1- c_{3}^2$; then $c_{i} \rightarrow 0,\quad  i = 1,2,3.$ Introduce spherical coordinates as follows: 
\begin{align*}
c_{1} = r \sin \varphi_1 \cos \varphi_2, \quad c_{2} = r \sin \varphi_1 \sin \varphi_2, \quad c_{3} = r \cos \varphi_1. 
\end{align*}
Then it follows from $r \rightarrow 0$ that:
\begin{align*}
&d_{i}^2 = 1 - (1-c_{3}^2)(1-c_{i}^2) = c_{i}^2 + c_{3}^2 - c_{i}^2 c_{3}^2 \newtilde c_{i}^2 + c_{3}^2, \quad i=1,2, \\
&\Delta = 1 - (1-c_{3}^2)(1-c_{1}^2)(1-c_{2}^2) \newtilde c_{1}^2 + c_{2}^2 + c_{3}^2 = r^2.
\end{align*}
This implies that 
\begin{align*}
\Big(\frac{d_{1} d_{2}}{\Delta}\Big)^2 &\newtilde \frac{(c_{1}^2 + c_{3}^2) (c_{2}^2 + c_{3}^2)}{r^4} = (\sin^2 \varphi_1 \cos^2 \varphi_2 + \cos^2 \varphi_1) (\sin^2 \varphi_1 \sin^2 \varphi_2 + \cos^2 \varphi_1) \leq 1,  \\
\Big(\frac{c_{1} d_{1}}{\Delta}\Big) &\newtilde \sin^2 \varphi_1 \cos^2 \varphi_2 (\sin^2 \varphi_1 \cos^2 \varphi_2 + \cos^2 \varphi_1) \leq 1.
\end{align*}
Therefore $\ds \frac{d_{1} d_{2}}{\Delta}, \frac{c_{1} d_{1}}{\Delta}$ are bounded from above.
\end{proof}

In the next lemmas we use the following parametrization of exponential mapping for the case $\lambda \in C_1$ (see~(\ref{ExpC1})):
\begin{align*}
x&=-\frac{4 \sigma k s_{1} s_{2} d_{1} d_{2}}{\alpha \Delta}, \\
y&=-\frac{4 \sigma}{\alpha} \bigg(\frac{k^2 s_{1} s_{2}^2 c_{1} d_{1}}{\Delta} + \frac{F_1}{2} - E_1\bigg), \\
z&=\frac{4 k c_{2} f_z}{|\alpha| \Delta}, \qquad \qquad \qquad \qquad f_z=c_{1} \big(F_1-2 E_1\big)+s_{1} d_{1}, \\
v&=\frac{y^3}{6}-\frac{2 k^2 (c_{1} c_{2} + s_{1} s_{2} d_{1} d_{2}) y}{|\alpha| \Delta^2} +\frac{4}{3\alpha \sigma} \bigg(F_1 \big(1-k^2\big) - E_1 \big(1- 2 k^2\big) \\
&-\frac{k^2 s_{1} d_{1}}{\Delta^3} \Big(6 s_{1} s_{2} c_{2} d_{1} d_{2} \big(2 d_{2}-\Delta \big) + c_{1} \big(1+3 c_{2}^2 (d_{2}^2-s_{2}^2)- k^4 s_{1}^2 s_{2}^6 (2 d_{1}^2 + s_{1}^2)\big)\Big)\bigg).
\end{align*}

\begin{lmm} \label{N11}
If $\{\nu_n\} \subset D_1 \cap N_1$ satisfies $\nu_n \rightarrow \bar{\nu} \in \cl(D_1) \backslash D_1$ ~ and ~ $\Exp (\nu_n) \rightarrow \bar{q} \in M_1$, then $\Delta \rightarrow 0$ or $\sigma \rightarrow 0$. 
\end{lmm}
\begin{proof}
Notice that for $\bar{\nu}=(\bar{x},\bar{y},\bar{z},\bar{v})\in M_1$ we have $\bar{x} \neq 0$ and $\bar{z} \neq 0$.
Consider all possible cases $\nu_n \rightarrow \partial D_1$:

\medskip 1. $\ds  \theta_{t/2} \rightarrow 0 ~ \Rightarrow ~ 
           \left\{  
           \begin{array}{rcl}  
           \sin \frac{\theta_{t/2}}{2} \rightarrow 0\\  
            \cos \frac{\theta_{t/2}}{2} \rightarrow 1
           \end{array}   
           \right. ~ \Rightarrow ~
           \left\{  
           \begin{array}{rcl}  
            k s_{2} \rightarrow 0, \\  
            d_{2} \rightarrow 1.   
           \end{array}   
           \right.$      
      It follows that $\ds x \newtilde \frac{(k s_{2}) s_{1} d_{1}}{\sigma} ~ \Rightarrow ~ \left[  
           \begin{array}{rcl}  
          	\sigma \rightarrow 0,\\  
            x \rightarrow \bar{x} = 0.   
           \end{array}   
           \right.
       $  
       
\medskip 2. $\ds \theta_{t/2} \rightarrow \pi ~ \Rightarrow ~  
           \left\{  
           \begin{array}{rcl}  
           \sin \frac{\theta_{t/2}}{2} \rightarrow 1\\  
            \cos \frac{\theta_{t/2}}{2} \rightarrow 0  
           \end{array}   
           \right. ~ \Rightarrow   ~ \left\{  
           \begin{array}{rcl}  
            k s_{2} \rightarrow 1, \\  
            d_{2} \rightarrow 0.   
           \end{array}   
           \right.$            
           Then $\ds x \newtilde \frac{s_{1} d_{1} d_{2}}{\sigma \Delta} ~ \Rightarrow ~ \left[  
           \begin{array}{rcl}  
          	\sigma \Delta \rightarrow 0,\\  
            x \rightarrow \bar{x} = 0.   
           \end{array}   
           \right.
           $      
            
\medskip 3. $\ds c_{t/2} \rightarrow 0 ~ \Rightarrow ~ k \sigma c_{2} \rightarrow 0 ~ \Rightarrow ~ ~ \ds z \newtilde \frac{k c_{2} f_z}{\sigma^2 \Delta} ~ \Rightarrow ~  \frac{\sigma^3 \Delta}{f_z} \rightarrow 0$, otherwise $z \rightarrow \bar{z} = 0$. 
This means that $\sigma^3 \Delta \rightarrow 0$ or $f_z \rightarrow \infty$.   
Suppose $f_z \rightarrow \infty, \sigma \rightarrow \bar{\sigma} \neq 0$, then  $u_1 \rightarrow \pi/2$ and $k \rightarrow 1$.                  
Since $k \sigma c_{2} \rightarrow 0$, then $u_2 \rightarrow \pi/2 ~ \Rightarrow ~ \Delta \rightarrow 0$.           

\medskip 4. $\ds  c_{t/2} \rightarrow \infty ~ \Rightarrow ~ k \sigma c_{2} \rightarrow \infty ~ \Rightarrow ~ \sigma \rightarrow \infty ~ \Rightarrow ~ \Delta \rightarrow 0,$ otherwise $x \rightarrow \bar{x} =  0.$

\medskip 5. $\ds t \rightarrow 0 ~ \Rightarrow ~ \frac{p}{\sigma} \rightarrow 0 ~ \Rightarrow ~ 
          \left[  
           \begin{array}{rcl}  
          	u_1 \rightarrow 0\\  
            \sigma \rightarrow \infty   
           \end{array}   
           \right. ~ \Rightarrow ~ \Delta \rightarrow 0,$ otherwise $x \rightarrow \bar{x} = 0.
$

\medskip 6. $\ds t \rightarrow \tmax ~ \Rightarrow 
          \left[  
           \begin{array}{rcl}  
          	f_z (u_1,k) \rightarrow 0 \text{ for } k \geq k_0\\  
            u_1 \rightarrow \pi \text{ for } k \leq k_0  
           \end{array}   
           \right. ~ \Rightarrow 
           \left[  
           \begin{array}{rcl}  
          	\sigma^2 \Delta \rightarrow 0, \text{ otherwise } z \rightarrow \bar{z} = 0 \\  
            \sigma \Delta \rightarrow 0, \text{ otherwise } x \rightarrow \bar{x} = 0.  
           \end{array}   
           \right.$
           
\medskip 7. $\ds |\alpha| \rightarrow \infty ~ \Rightarrow ~ \sigma \rightarrow \infty ~ \Rightarrow ~ \Delta \rightarrow 0,$ otherwise $x \rightarrow \bar{x} = 0$.           
\end{proof}

\begin{lmm} \label{N12}
Suppose $\nu_n \in D_1 \cap N_{1}$. If $\Delta \rightarrow 0$, then $x\rightarrow 0$ or $y \rightarrow \infty$.
\end{lmm}
\begin{proof}
Consider two possible cases:

\medskip 1. $\sigma \rightarrow \infty ~ \Rightarrow ~ \ds x = - 4 k s_{1} s_{2} \frac{d_{1} d_{2}}{\Delta} \frac{1}{\sigma} \rightarrow 0$ (see~Lemma~\ref{lemm1}).

\medskip 2. $\sigma \rightarrow \bar{\sigma} < \infty$. It follows from Lemma~\ref{lemm1} that $\ds k^2 s_{1} s_{2}^2 \frac{c_{1} d_{1}}{\Delta} - E_1$ is bounded from above. Since $F_1 \rightarrow \infty$, then $y \rightarrow \infty$.
\end{proof}

\begin{lmm} \label{N13}
Suppose $\nu_n \in D_1 \cap N_{1}$. If $\sigma \rightarrow 0, \Delta \rightarrow \bar{\Delta} \neq 0$, then one of the functions $x,y,z,v$ tends to $\infty$, otherwise $x$ or $z$ tends to zero.
\end{lmm}
\begin{proof}
Assume the converse. Then notice that $\ds k s_{1} s_{2} d_{1} d_{2} \rightarrow 0$, otherwise $\ds x \newtilde \frac{1}{\sigma} \rightarrow \infty$. The proof consists of the following six items:

\medskip 1. $\ds d_{1} \rightarrow 0$. This means that $\ds u_1 \rightarrow \pi/2, k\rightarrow 1$ and $\ds u_2 \rightarrow \overline{u_2} \neq \pi/2$. Whence, $\ds y \newtilde F_1/\sigma \rightarrow \infty$.

\medskip 2. $u_1 \rightarrow 0$. Consider four subcases:

\medskip 2.1.  $s_{2} k \rightarrow 0$. Here we have  
\begin{align*}
&x \newtilde \frac{s_{2} k u_1}{\sigma} ~ \Rightarrow ~ s_{2} k u_1 \newtilde \sigma, \\
&F_1 \sim u_1, E_1 \sim u_1 ~ \Rightarrow ~ y \newtilde \frac{1}{\sigma} \big(s_{1} (s_{2} k)^2 \frac{c_{1} d_{1}}{\Delta} + \frac{F_1}{2} - E_1 \big) \newtilde \frac{u_1}{\sigma} \newtilde \frac{1}{s_{2} k} \rightarrow \infty.
\end{align*}

\medskip 2.2. $d_{2} \rightarrow 0$. It follows that
\begin{align*}
&x \newtilde \frac{d_{2} u_1}{\sigma} ~ \Rightarrow ~ d_{2} u_1 \newtilde \sigma, \qquad &y \newtilde \frac{u_1}{\sigma} \newtilde \frac{1}{d_{2}} \rightarrow \infty.
\end{align*}

\medskip 2.3. $c_{2} \rightarrow 0, \ k \rightarrow \bar{k} \in (0, 1)$. We get
\begin{align*}
&x \newtilde \frac{u_1}{\sigma} ~ \Rightarrow ~ u_1 \newtilde \sigma, \qquad &z \newtilde \frac{c_{2} u_1^3}{\sigma^2} \newtilde c_{2} u_1 \rightarrow 0.
\end{align*}

\medskip 2.4. $k \rightarrow \bar{k} \neq 0, u_2 \rightarrow \overline{u_2} \in (0, \pi/2)$. Hence
\begin{align*}
&x \newtilde \frac{u_1}{\sigma}, \quad &z \newtilde \frac{u_1^3}{\sigma^2} \newtilde x^2 u_1 \rightarrow 0.
\end{align*}

\medskip 3. $u_1 \rightarrow \pi ~ \Rightarrow ~ F_1 \rightarrow 2 K, \ E_1 \rightarrow 2 E$. We obtain
$\ds y \newtilde \frac{2 E - K}{\sigma} ~ \Rightarrow ~ 2 E - K \rightarrow 0$. It follows that $\ds v \newtilde \frac{K}{\sigma^3} \rightarrow \infty$, since $K \rightarrow \overline{K} = K(\bar{k}) > 0$.

\medskip 4. $\ds d_{2} \rightarrow 0, \ u_1 \rightarrow \overline{u_1} \in (0, \pi/2) ~ \Rightarrow ~ u_2 \rightarrow \pi/2, \ k \rightarrow 1$. 
\begin{align*}
&x \newtilde \frac{d_{2}}{\sigma}, \qquad &y \newtilde \frac{1}{\sigma} \Big(s_{1} + \frac{F_1}{2} - E_1 \Big).
\end{align*}
Note that the function $\ds \sin u_1 + \frac{F(u_1,1)}{2} - E(u_1,1)$ vanishes only at the point $u_1 = 0$ since it has positive derivative $\ds \frac{1}{2 \sqrt{1-\sin^2 u_1}}$. Therefore $y \rightarrow \infty$.

\medskip 5. $\ds u_2 \rightarrow 0, \ s_{1} \rightarrow \overline{s_{1}} \neq 0, \ d_{1} \rightarrow \overline{d_{1}} \neq 0, \ k \rightarrow \bar{k} \neq 0$. We have
$y \newtilde \frac{2 E_1 - F_1}{\sigma} ~ \Rightarrow ~ 2 E_1 - F_1 \rightarrow 0$. Hence $\ds z \newtilde \frac{1}{\sigma^2} \rightarrow \infty$.

\medskip 6. $k \rightarrow 0$.  We get $\ds y \newtilde \frac{u_1}{\sigma} ~ \Rightarrow ~ u_1 \rightarrow 0$ (see~item~2.1).
\end{proof}

In the next lemmas we use the following parametrization of exponential mapping for the case $\lambda \in C_2$ (see~(\ref{ExpC2})):
\begin{align*}
x&=-\frac{4 \sgn c\ \sigma s_{1} s_{2} c_{1} c_{2}}{\alpha k \Delta}, \\
y&=-\frac{4 \sigma}{\alpha k} \bigg(\frac{k^2 s_{1} s_{2}^2 c_{1} d_{1}}{\Delta} + \Big(1- \frac{k^2}{2}\Big)F_1 - E_1\bigg), \\
z&=-4\frac{4 \sgn c\ d_{2} g_z} {|\alpha| k^2 \Delta}, \qquad g_z=\big(2 E_1 +(k^2-2)F_1\big) d_{1} - k^2 s_{1} c_{1}, \\
v&=\frac{y^3}{6}+ \frac{2 y}{|\alpha| k^2 \Delta^2}\bigg(1 + (1 + c_{1}^2 c_{2}^2) k^4 s_{1}^2 s_{2}^2 - k^2 (s_{1}^2 + s_{2}^2 - 2 c_{1} c_{2} d_{1} d_{2} s_{1} s_{2})\bigg)\\  
&  - \frac{3}{4 \alpha \sigma k}\bigg( 2 F_1\Big(\frac{1}{k^2}-1\Big) - E_1 \Big(\frac{2}{k^2} - 1\Big) + \frac{c_{1} s_{1}}{\Delta^3} \Big(2 c_{2}^2 d_{1} (1+d_{1}^2) \Delta + 6 c_{1} c_{2} d_{2} k^2 s_{1} s_{2} (2 c_{2}^2 - \Delta)  \\
&+ d_{1} s_{2}^2 \big((2-k^2)\Delta^2- 4 d_{1}^2 d_{2}^2\big)\Big)\bigg).             
\end{align*}

\begin{lmm} \label{N21}
If $\{\nu_n\} \subset D_1 \cap N_{2}$ satisfies $\nu_n \rightarrow \bar{\nu} \in \cl(D_1) \backslash D_1$ and $\Exp (\nu_n) \rightarrow \bar{q} \in M_1$, then $\Delta \rightarrow 0$ or $\ds \frac{\sigma}{k} \rightarrow 0$. 
\end{lmm}
\begin{proof}
Notice that for $\bar{\nu}=(\bar{x},\bar{y},\bar{z},\bar{v})\in M_1$ we have $\bar{x} \neq 0$. Consider all possible cases for $\bar{\nu} \in \cl(D_1) \backslash D_1$:

\medskip 1. $\ds \theta_{t/2} \rightarrow 0 ~ \Rightarrow ~  
           \left\{  
           \begin{array}{rcl}  
           \sin \frac{\theta_{t/2}}{2} \rightarrow 0\\  
            \cos \frac{\theta_{t/2}}{2} \rightarrow 1
           \end{array}   
           \right. ~ \Rightarrow ~
           \left\{  
           \begin{array}{rcl}  
            s_{u_2} \rightarrow 0, \\  
            c_{u_2} \rightarrow 1.   
           \end{array}   
           \right.$           
         It follows that $\ds \frac{\sigma \Delta}{k} \rightarrow 0$, otherwise $x \rightarrow \bar{x} = 0$.

\medskip 2. $ \ds \theta_{t/2} \rightarrow \pi ~ \Rightarrow ~ 
           \left\{  
           \begin{array}{rcl}  
           \sin \frac{\theta_{t/2}}{2} \rightarrow 1\\  
            \cos \frac{\theta_{t/2}}{2} \rightarrow 0  
           \end{array}   
           \right. ~ \Rightarrow ~
           \left\{  
           \begin{array}{rcl}  
            s_{2} \rightarrow 1, \\  
            c_{2} \rightarrow 0.   
           \end{array}   
           \right.$
           This means that $\ds \frac{\sigma \Delta}{k} \rightarrow 0$, otherwise $ x \rightarrow \bar{x} = 0$.       

\medskip 3. $ \ds  c_{t/2} \rightarrow 0 ~ \Rightarrow ~ \frac{\sigma}{k} d_{2} \rightarrow 0 ~ \Rightarrow ~ \frac{\sigma}{k} \rightarrow 0$ or $c_{2} \rightarrow 0$.
				From  $\ds c_{2} \rightarrow 0$ we have $\ds \frac{\sigma \Delta}{k} \rightarrow 0$, otherwise $x \rightarrow \bar{x} = 0$.

\medskip 4. $ \ds c_{t/2} \rightarrow \infty ~ \Rightarrow ~ \frac{\sigma}{k} d_{2} \rightarrow \infty ~ \Rightarrow ~ \frac{\sigma}{k} \rightarrow \infty ~ \Rightarrow ~ \Delta \rightarrow 0$, otherwise $x \rightarrow \bar{x} = 0$.

\medskip 5. $ \ds t \rightarrow 0 ~ \Rightarrow ~ \frac{p k}{\sigma} \rightarrow 0 ~ \Rightarrow ~ \frac{s_{1} k}{\sigma} \rightarrow 0 ~ \Rightarrow ~ \Delta \rightarrow 0$, otherwise $x \rightarrow \bar{x} = 0$.

\medskip 6. $ \ds t \rightarrow \tmax ~ \Rightarrow ~ u_1 \rightarrow \frac{\pi}{2} ~ \Rightarrow ~ \frac{\sigma \Delta}{k} \rightarrow 0$, otherwise $x \rightarrow \bar{x} = 0$.

\medskip 7. $ \ds \alpha \rightarrow \infty ~ \Rightarrow ~ \sigma \rightarrow \infty ~ \Rightarrow ~ \Delta \rightarrow 0$, otherwise $x \rightarrow \bar{x} = 0$.
\end{proof}

\begin{lmm} \label{N22}
Suppose $\nu_n \in D_1 \cap N_{2}$. If $\Delta \rightarrow 0$, then $x\rightarrow 0$ or $y \rightarrow \infty$.
\end{lmm}
\begin{proof}
Consider two possible cases:

\medskip 1. $\sigma \rightarrow \infty ~ \Rightarrow ~ \ds x = - 4 k s_{1} s_{2} \frac{c_{1} c_{2}}{\Delta} \frac{1}{\sigma} \rightarrow 0$ (see~Lemma~\ref{lemm1}).

\medskip 2. $\sigma \rightarrow \bar{\sigma} < \infty$. It follows from~Lemma~\ref{lemm1} that $\ds\big(k^2 s_{1} s_{2}^2 \frac{c_{1} d_{1}}{\Delta} - E_1\big)$ is bounded from above. And since $F_1 \rightarrow \infty$ we have $y \rightarrow \infty$.
\end{proof}

\begin{lmm} \label{N23}
Suppose $\nu_n \in D_1 \cap N_{2}$. If $\frac{\sigma}{k} \rightarrow 0, \Delta \rightarrow \bar{\Delta} \neq 0$, then one of the functions $x,y,z$ or $v$ tends to $\infty$, otherwise $x$ or $z$ tends to 0.
\end{lmm}
\begin{proof}
Assume the converse. Then notice that $s_{1} s_{2} c_{1} c_{2} \rightarrow 0$, otherwise $\ds x \newtilde \frac{k}{\sigma} \rightarrow \infty$. The proof consists of five steps:

\medskip 1. $u_1 \rightarrow 0$. Then we obtain 
\begin{align*}
&x \newtilde \frac{u_1 s_{1} c_{2}}{\sigma/k} ~ \Rightarrow ~ \frac{\sigma}{k} \newtilde u_1 s_{2} c_{2}. 
\end{align*}  
It follows from Taylor expansion that $g_z \newtilde k^2 u_1^3$, then
\begin{align*}
&z \newtilde \frac{d_{2} k^2 u_1^3}{\sigma^2} \newtilde \frac{d_{2} u_1}{s_{2}^2 c_{2}^2} = \frac{d_{2} u_1}{s_{2}^2 c_{2}^2} ~ \Rightarrow ~ s_{2} c_{2} \rightarrow 0,
\end{align*}  
otherwise $z \rightarrow 0$. 
\begin{align*}
&y \newtilde \frac{1}{\sigma k} \Big(k^2 u_1 s_{2}^2 + (1- k^2/2)u_1-u_1\Big) = \frac{u_1 (s_{2}^2 - 1/2)}{\sigma/k} \newtilde \frac{s_{2}^2 - 1/2}{s_{2} c_{2}} \rightarrow \infty.
\end{align*}  

\medskip 2. $k\rightarrow 0, \ u_1 \rightarrow \overline{u_1} \neq 0$. Using Taylor expansion we get $\ds z \newtilde \frac{k^2}{\sigma^2} \rightarrow \infty$.

\medskip 3. $\ds u_1 \rightarrow \frac{\pi}{2}, k \rightarrow \bar{k} \neq 0$. We have $\ds y \newtilde \frac{1}{\sigma} \Big(\frac{k^2 c_{1} s_{1} d_{1}^2 s_{2}^2}{\Delta} + (1-k^2/2)F_1 - E_1 \Big).$
Notice that $$\ds \frac{d}{d k}\Big((1-k^2/2)F_1 - E_1\Big) = \frac{k}{2(1-k^2)} \int_0^{\pi/2} \frac{k^2 \cos^2 \theta d\theta}{\sqrt{1-k^2 \sin^2 \theta}}> 0 ~ \Rightarrow ~ (1-k^2/2) F_1 - E_1 > 0.$$
Combining the last inequality and $\ds \frac{k^2 c_{1} s_{1} d_{1}^2 s_{2}^2}{\Delta} \rightarrow 0$, we obtain $y \rightarrow \infty$.

\medskip 4. $u_2 \rightarrow 0, \ u_1 \rightarrow \overline{u_1} \in (0, \pi/2), \ k \rightarrow \bar{k} \neq 0$. Here we have
\begin{align*}
&x \newtilde \frac{u_2}{\sigma} \Rightarrow u_2 \newtilde \sigma, 
\end{align*}
otherwise $x \rightarrow 0$. 
\begin{align*}
&z \newtilde \frac{(2 E_1 + (k^2 -2) F_1) d_{1} - k^2 c_{1} s_{1}}{\sigma^2}, \\
&y \newtilde \frac{1}{\sigma} \Big( k^2 s_{1} c_{1} d_{1} s_{2}^2 + (1-k^2/2) F_1 - E_1 \Big).
\end{align*}
Since $k^2 s_{1} c_{1} d_{1} s_{2}^2 \rightarrow 0$, we see that $(1-k^2/2) F_1 - E_1 \rightarrow 0$, otherwise $y\rightarrow \infty$. Hence from $k^2 c_1 s_1 \newtilde 1$ we get
$\ds z \newtilde \frac{1}{\sigma^2} \rightarrow \infty$.

\medskip 5. $u_2 \rightarrow \pi/2, u_1 \rightarrow u_1 \in (0, \pi/2), k \rightarrow \bar{k} \neq 0$. Suppose $\bar{k} \neq 1$, then $\ds z\newtilde \frac{d_{2} g_z}{\sigma^2} \rightarrow \infty$. This means that $\bar{k} = 1$. Here we have $\ds y \newtilde \frac{1}{\sigma} \Big( \frac{c_{1} s_{1} d_{1}}{\Delta} + F_1/2 - E_1 \Big) = \frac{1}{\sigma} (s_{1} + F_1/2 - E_1)$.
Since $\ds \frac{d}{d u_1} (s_{1} + F_1/2 - E_1) = \frac{1}{2 \cos u_1} > 0$, it follows that $y \rightarrow \infty$.
\end{proof}

\begin{thrm}\label{proper}
The mapping $\Exp : D_i \rightarrow M_i$ is proper for $i = 1,\ldots,4$.
\end{thrm}
\begin{proof}
Assume the converse. Then it follows from Lemma~\ref{lm1} that $\Exp : D_1 \rightarrow M_1$ is not proper. By Lemma~\ref{lm2}, there exists a sequence $\nu_n \in \newD$, such that $\nu_n \rightarrow \bar{\nu} \in \cl (D_1) \backslash D_1, ~ \Exp(\nu_n) \rightarrow \bar{q} \in M_1$. Since  $\nu_n \in \newD$, we consider 3 cases:
\begin{enumerate}
\item $\{\nu_n\} \subset D_1 \cap N_{6}$   is impossible (see~Lemmas~\ref{N61},~\ref{N62}), 
\item $\{\nu_n\} \subset D_1 \cap N_{3}$ is impossible (see~Lemmas~\ref{N31}--\ref{N33}),
\item $\{\nu_n\} \subset (D_1 \cap N_{1}) \cup (D_1 \cap N_{2})$ is impossible (see~Lemmas~\ref{N11}--\ref{N23}).
\end{enumerate}
Since all cases are impossible, we have a contradiction which proves the theorem.
\end{proof}

\begin{thrm} \label{Expdif}
The mapping $\Exp: D_i \rightarrow M_i $ is a diffeomorphism for $i=1,\ldots,4$. 
\end{thrm}
\begin{proof}
Follows from Th.~\ref{Hadamard}, since all hypotheses of this theorem hold by Propos.~\ref{DiMitopo}, Th.~\ref{th:tconjmax} and Th.~\ref{proper}.
\end{proof}
\begin{crllr}\label{Expdiffeo}
The mapping $\Exp: \widetilde{N} \rightarrow \widetilde{M}$ is a diffeomorphism.
\end{crllr}

\section{Cut time} 
In this section we prove that the cut time coincides with the first Maxwell time corresponding to reflections.
\subsection{Cut time and Maxwell time}
\begin{thrm} \label{th:tcut}
For any $\lambda \in C$, $$\tcut (\lambda) = \tmax (\lambda).$$
\end{thrm}

\begin{proof}
Take any $\lambda \in C$ and denote $t_1 = \tmax (\lambda)$. Since $\tcut (\lambda) \leq t_1$ by Th.~\ref{th:tcut_bound}, it remains to prove that $\tcut (\lambda) \geq t_1$.

Let us call a pair $(\lambda , t) \in N$ optimal if the geodesic $\Exp(\lambda, s)$ is optimal on the segment $s \in [0, t]$. We have to show that $(\lambda, t)$ is optimal for any $t \in (0, t_1)$.

\begin{enumerate}

\item If $\lambda \in C_4 \cup C_5 \cup C_7$, then $t_1 = +\infty$, and any $(\lambda, t), \ t \in (0, t_1)$, is optimal since $(x_s,y_s)$ is a straight line.
\item\label{itt2} Let $\lambda \in C_1 \cup C_2 \cup C_6 ,$ thus $t_1 \in (0, +\infty)$.

Since $t_1 = \tmax (\lambda)$, then $ \nu_1 = (\lambda, t_1) \in N'$. For $\lambda \in C_1 \cup C_2 \cup C_6$ the function $t \mapsto \sin {\theta}_{t/2} c_{t/2}$ has isolated zeros, thus there exists $t \in (0, t_1)$ arbitrarily close to $t_1$ such that $\nu = (\lambda, t) \in \widetilde{N}$. Then $q = \Exp (\nu) \in \widetilde{M}$ (Propos.~\ref{ExpDiMi}). Since $\Exp(N') \cap \widetilde{M} = \emptyset$ (Propos.~\ref{ExpN'}) and $\Exp \colon \widetilde{N} \to \widetilde{M}$ is a diffeomorphism (see Corollary~\ref{Expdiffeo}) then $\Exp^{-1} (q) \cap \widehat{N} = \{\nu\}$. Thus $\nu = (\lambda, t)$ is optimal. Since $t$ can be chosen arbitrarily close to $t_1$, then any $(\lambda, t), t \in (0, t_1)$, is optimal.
\item Let $\lambda \in C_3$, then $t_1 = + \infty$. There exist $(\lambda,t) \in \widetilde{N}$ for arbitrarily large $t$. Then the proof follows the argument of item~(\ref{itt2}).
\end{enumerate}
\end{proof}

Now we collect all properties of the cut time that we previously obtained for the Maxwell time $\tmax$.


\begin{crllr}\label{tcut} The function $\tcut \colon C \to (0, +\infty]$ has the following properties:

\begin{enumerate}
\item Let $\lambda \in C$ and let $t_1 = \tcut (\lambda)$. For finite $t_1$, a trajectory $\Exp(\lambda, s), \ s \in [0, t],$ is optimal iff $t \in [0, t_1]$. For $t_1 = +\infty$, any trajectory $\Exp(\lambda, s), \ s \in [0, t], \ t > 0$, is optimal. 

\item The function $\tcut$ has the following explicit representation:
\begin{align*}
& \forall \lambda \in C_1 &\quad & \tcut (\lambda) = \frac{p_1 (k)}{ \sqrt{ |\alpha| } }, \\
& \forall \lambda \in C_2 &\quad & \tcut (\lambda) = \frac{2 K k}{ \sqrt{ |\alpha| } }, \\
& \forall \lambda \in C_6 &\quad & \tcut (\lambda) = \frac{2 \pi}{ \sqrt{|c| } }, \\
& \forall \lambda \in C_3 \cup C_4 \cup C_5 \cup C_7 &\quad & \tcut (\lambda) = +\infty. 
\end{align*}

\item The function $\tcut$ depends only on $E$ and $|\alpha|$, is preserved by the flow of $\vec{H}_v$ and by the reflections ${\varepsilon}^i$, and is homogeneous of order one w.r.t. the dilations $\delta_\mu$. 

\item The function $\tcut$ is continuous on $C \backslash C_4$ and is smooth on $C^0_1 \cup C_2$, where $C^0_1 = \{ \lambda \in C_1 ~ | ~ k \ne k_0 \}$.

\end{enumerate}
\end{crllr}


According to Cor.~\ref{tcut}, the function $\tcut  \colon C \to (0, +\infty ]$ is invariant w.r.t. the flow $e^{s \vec{H}_v}$ and the reflections ${\varepsilon}^i$, and respects the action of dilations $\delta_\mu$:
\begin{align*}
& \tcut \circ e^{s \vec{H}_v} = \tcut \circ {\varepsilon}^i = \tcut, \\
& \tcut \circ \delta_\mu = \mu \tcut. 
\end{align*}
Thus the cut time can be represented (up to a constant positive factor) by a univariate function on the quotient $ C / \left\langle e^{s \vec{H}_v}, {\varepsilon}^i, \delta_\mu \right\rangle$. The quotient $C / \left\langle e^{s \vec{H}_v} , {\varepsilon}^i \right\rangle$ can be represented by the quadrant $\{ (\theta , c, \alpha) \in C ~ | ~ \theta = 0, \ c \geq 0, \ \alpha \geq 0 \}, $ thus 
$$ C / \left\langle  e^{s \vec{H_v}} , {\varepsilon}^i , \delta_\mu \right\rangle \cong \Gamma \sqcup P , $$
 where 
\begin{align*}
& \Gamma = \{ (\theta , c, \alpha) \in C ~ | ~ \theta = 0, \ c = \sin \beta, \ \alpha = \cos \beta, \ \beta \in [0, \pi/2] \}, \\
& P = \{ (\theta , c, \alpha) \in C ~ | ~ \theta = 0, \  c = 0, \ \alpha = 0 \}. 
\end{align*}

The point $P$ corresponds to the subset $C_7$,
while the arc $\Gamma$ corresponds to the rest subsets $C \backslash C_7 = C_1 \cup C_2 \cup C_{35} \cup C_4 \cup C_6$ of decomposition~(\ref{C1357}).

Thus (up to a constant positive factor) the cut time can be represented on the set $C \backslash C_7$ as a univariate function $ \tcut (\beta), \ \beta \in [0, \pi/2]$.

If $\beta = 0,$ then $\lambda \in C_4$, thus $\tcut (\beta) = +\infty$.

If $\beta \in (0, {\beta}_1)$, where ${\beta}_1 = \arccos ( \sqrt{5} - 2)$, then $\lambda \in C_1$, thus
\begin{align*}
& \tcut (\beta) = \frac{2 p_1 (k)}{ \sqrt{\alpha} }, \quad k = \sqrt{ \frac{ {\sin}^2 \beta }{ 4 \cos \beta } }, \quad \alpha = \cos \beta, \quad p_1 (k) = \min \big( p_z^1 (k), 2 K (k) \big).
\end{align*} 

If $\beta = {\beta}_1$, \ie, ${\sin}^2 \beta = 4 \cos \beta$, then $\lambda \in C_{35}$, thus $\tcut (\beta) = +\infty$.

If $\beta \in ({\beta}_1 , \pi/2)$, then $\lambda \in C_2$, thus 

$$ \tcut (\beta) = \frac{2 K k}{ \sqrt{\alpha} } , \quad k = \sqrt{ \frac{ 4 \cos \beta }{ {\sin}^2 \beta }}, \quad \alpha = \cos \beta. $$


Finally, if $\beta = \pi/2$, then $\lambda \in C_6$, thus $\tcut (\beta) = 2 \pi$.

\begin{figure}[htbp]
\begin{center}
\includegraphics[width=0.75\linewidth]{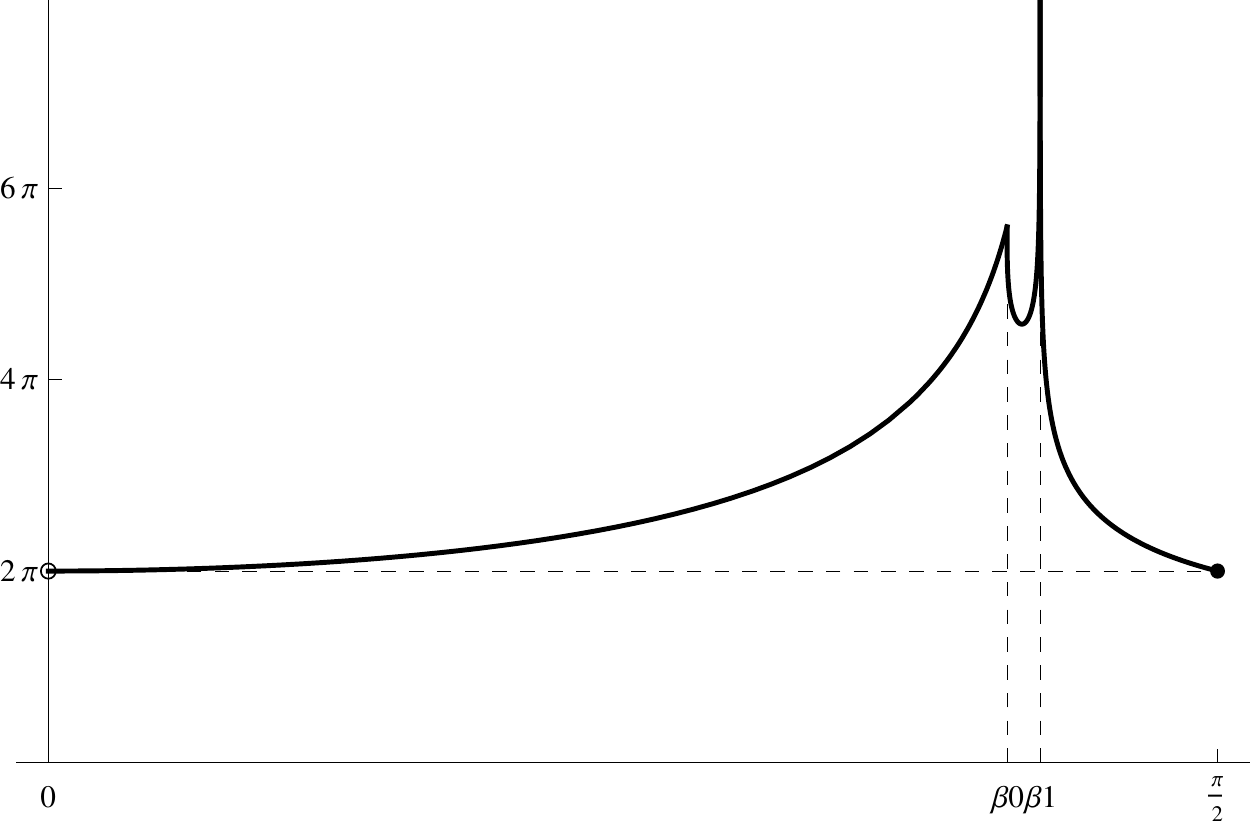}
\parbox{0.9\linewidth}{\caption{Plot of the function $\beta \mapsto \tcut(\beta)$}\label{tcutbeta}}
\end{center}
\end{figure}

The plot of the function $\tcut (\beta)$ is shown in Fig.~\ref{tcutbeta}.
Notice continuity of $\tcut  (\beta)$ everywhere except $\beta = 0$, where $\tcut (+0) = 2 \pi < +\infty = \tcut (0)$.
Also notice smoothness of $\tcut (\beta)$ everywhere except $\beta = 0, \ \beta = {\beta}_1$, and $\beta = {\beta}_0$, where $\ds \frac{ {\sin}^2 {\beta}_0 }{ 4 \cos {\beta}_0 } = k^2_0 , \ 2 E (k_0) - K (k_0) = 0$ (here $E (k_0)$ is the complete elliptic integral of the second kind), corresponds to the figure-of-eight closed Euler elastica. These regularity properties of $\tcut (\lambda) = \tmax (\lambda)$ are reported in Corollary~\ref{tcut}.


\subsection{Cut time and conjugate time}

\begin{prpstn} \label{tcutconjC1}
Let $\lambda \in C_1, \ t_1 = \tcut (\lambda), \ \tau = ( \varphi + t_1 / 2 ) / \sqrt{ |\alpha| }$. Then $t_1 = \tconj (\lambda)$ iff one of the conditions hold:
\begin{enumerate}
\item $k < k_0, \ \sn \tau = 0,$
\item $k = k_0 ,$
\item $k > k_0, \ \cn \tau = 0.$
\end{enumerate}

In particular, if $t_1 = \tconj (\lambda)$, then $\sn \tau \cn \tau = 0$ or $k = k_0$.
\end{prpstn}

\begin{proof} 
Follows immediately from Lemma 8~\cite{engel_conj}.
\end{proof}

\begin{rmrk}
The equality $\sn \tau = 0 \ (\cn \tau = 0)$ is equivalent to $\sin {\theta}_{t/2} = 0 \ (c_{t_1/2} = 0);$ it means that elastica $(x_t , y_t), \ t \in [0, t_1]$, is centered at a vertex (resp. inflexion point). The equality $k = k_0$ means that $(x_t , y_t), \ t \in [0, t_1]$, is the closed figure-of-eight elastica.
\end{rmrk}

\begin{prpstn} \label{tcutconjC2}
Let $\lambda \in C_2, \ t_1 = \tcut (\lambda), \ \tau = ( \varphi + t_1 / 2 ) / (k \sqrt{|\alpha|})$. Then $t_1 = \tconj (\lambda)$ iff $\sn \tau \cn \tau = 0$. 
\end{prpstn}
\begin{proof}
Follows immediately from Lemma 8~\cite{engel_conj}.
\end{proof}

\begin{rmrk} 
The equality $\sn \tau \cn \tau = 0$ is equivalent to $\sin  {\theta}_{t/2}  = 0$. It means that the corresponding elastica is centered at vertex.
\end{rmrk}


\begin{prpstn} \label{tucutconjC6}
Let $\lambda \in C_6 , \quad t_1 = \tcut (\lambda)$. Then $t_1 = \tconj (\lambda)$ iff $\sin \theta = 0.$
\end{prpstn}

\begin{proof}
Let $\lambda = (\theta , c, \alpha) \in C_6 , \ \alpha = 0, \ c \ne 0, t_1 = \tcut (\lambda), \nu = (\lambda, t_1) \in N_6$. Since $C_6 \subset \cl (C_2)$, the expression for Jacobian  $ \ds J = \frac{ \partial (x, y, z, v) }{ \partial (\theta, c, \alpha, t) } $ for $\nu \in N_6$
 can be obtained by passing to the limit $\alpha \to 0$ in the expression for Jacobian $J|_{_{N_2}}$ computed in~\cite{engel_conj}. By such a limit we get $\ds J(\nu) = \frac{ {\pi}^3 }{ {|c|}^3 } {\sin}^2 \theta $. So the instant $t_1$ is a conjugate time iff $\sin \theta = 0$.
\end{proof} 
 
\subsection{Optimal trajectories for special boundary conditions}
 
For a generic terminal point $q_1 \in \widetilde{M}$, there exists a unique optimal trajectory $q_t = \Exp(\lambda, t), \ t \in [0, t_1]$, which can be found by solving the equation $\Exp(\lambda , t_1) = q_1 , \ (\lambda , t_1) \in \widetilde{N}$. 
 
 In this subsection we discuss special boundary conditions for which optimal trajectories can be given explicitly or by a more simple equation.
 
 \subsubsection{Abnormal variety}
 
 Consider the set of points in $M$ filled by abnormal trajectories:
 $$ A = \{ q \in M \mid x = z = 0, \ v=y^3 / 6 \}. $$
 We have $\Exp(C_4 , {\mathbb{R}}_+) = \Exp (C_5 , {\mathbb{R}}_+) = \Exp (C_7^{0, \pi}, {\mathbb{R}}_+ ) = A \backslash \{ q_0 \}$, where 
 $$C_7^{0, \pi} = \big\{ \lambda = (\theta, c, \alpha) \in C_7 \mid \alpha = c = 0, \ \theta \in \{0, \pi\} \big\}.$$ 
 Any nonzero point $q_1 = (0, y_1 , 0, v_1) \in A$ is connected with $q_0$ by a unique optimal trajectory
 $$ x_t = 0, \quad y_t = t \sgn y_1, \quad z_t = 0, \quad v_t = \frac{t^3}{6} \sgn y_1, \qquad t \in [0, |y_1|]. $$
 
 \subsubsection{Straight lines $(x_t , y_t)$}
 
 The set of points in $M$ filled by trajectories that project to straight lines $(x_t , y_t)$ is 
 
 $$ L = \{ q \in M \mid z=0 , v = (x^2 + y^2) y / 6 \} \supset A. $$
 We have $\Exp(C_7 , {\mathbb{R}}_+) = L \backslash \{ q_0 \}$.
 The unique optimal trajectory for $q_1 \in L \backslash \{ q_0 \}$ is $\Exp(\lambda , t), \lambda \in C_7 ,$ \ie, 
 $$ x_t = -t \sin \theta, \quad y_t = t \cos \theta, \quad z_t = 0, \quad v_t = \frac{t^3}{6} \cos \theta, \qquad t \in [0, t_1],$$ 
where $t_1 > 0$ and $\theta \in S^1$ are found from the equations 
 $$ x_1 = -t_1 \sin \theta, \qquad y_1 = t_1 \cos \theta.$$
 
 
 \subsubsection{Fixed points of reflection ${\varepsilon}^6$}
 
The reflection ${\varepsilon}^6 \colon M \to M$ has the set of fixed points
$$ S_6 = \{ q \in M \mid y=0, \ v=xz/2 \}.$$
This 2-dimensional manifold is of particular interest since it is the only fixed manifold of reflections
$$ S_i = \{ q \in M \mid {\varepsilon}^i (q) = q \}, \quad i = 1, \dots , 7, $$
not contained completely in codimension one manifolds $S_1 = \{ q \in M \mid z = 0 \} $ and $ S_2 = \{ q \in M \mid x = 0 \} $. By Lem.~4~\cite{engel}, 
$$ S_3 \cup S_4 \cup S_5 \subset S_2 , \qquad S_7 \subset S_1. $$
But $S_6 \not\subset S_1 \cup S_2$. 

If a point $q_1 \in S_6 \backslash \{ q_0 \}$ is connected with $q_0$ by a trajectory $q_t = \Exp (\lambda , t), \ t \in [0, t_1], \ q_{t_1} = q_1$, then the trajectory $q^6_t = \Exp({\lambda}^6 , t), \ t \in [0, t_1]$, satisfies the equation $q_{t_1}^6 = q_1$. 

Moreover, if ${\lambda}^6 \ne \lambda$, then the 
points $(\lambda, t_1), ({\lambda}^6 , t_1)$ would belong to the Maxwell set 
$$ {\MAX}^6 = \{ (\lambda, t) \in N \mid {\lambda}^i \ne \lambda, \ \Exp({\lambda}^i , t) = \Exp(\lambda, t) \}, $$
which might a priori give Maxwell times which are additional to those provided by the sets ${\MAX}^1 , {\MAX}^2$ studied in~\cite{engel}. It turns out that, as we show below, the equality ${\lambda}^6 = \lambda$ is satisfied for all $\lambda \in C$ with $\Exp (\lambda, t_1) \in S_6, \ t_1 > 0$.
 
Consider the decomposition
\begin{align*}
& S_6 = \bigsqcup_{i, j \in \{ 0, +, - \}} S_{ij} , \\
& S_{ij} = \{ q \in S_6 \mid \sgn x = i, \ \sgn z = j \}. 
\end{align*} 
 For example, $S_{+-} = \{ q \in S_6 \mid x > 0, z < 0 \}.$
 
 Denote $ N_{ij} = \{ (\lambda, t) \in N_6 \mid \tau = -i \pi /2, \ \sgn c = j, \ t \in (0, 2 \pi / |c| ) \}, \quad i, j \in \{ +, - \}, $
 where $\tau = \theta + c t/2$.
 
\begin{lmm} \label{ExpNij}
For any $i, j \in \{ +, i \}$, the mapping $\Exp \colon N_{ij} \to S_{ij}$
is a diffeomorphism.
\end{lmm}
 
\begin{proof} The reflections ${\varepsilon}^4$ and ${\varepsilon}^7$ permute the sets $N_{ij} , S_{ij}$, thus it suffices to prove only that the mapping $\Exp \colon N_{++} \to S_{++}$ is a diffeomorphism.
 
If $(\lambda, t) \in N_{++}$, then 
\begin{align}
x_t = \frac{2 \sin p}{ c}, \quad y_t=0, \quad z_t = \frac{2p - \sin (2p)}{2 c^2}, \quad v_t = \frac{x_t z_t}{2}, \label{ExpN++}
\end{align}
where $ p = c t /2 \in (0, \pi) $. Thus $\Exp(N_{++}) \subset S_{++}$. 

Further, the mapping $\Phi \colon (p, c) \mapsto (x,z)$ is a diffeomorphism from $(0, \pi) \times {\mathbb{R}}_+$ to ${\mathbb{R}}_+ \times {\mathbb{R}}_+$ by Hadamard global diffeomorphism theorem, thus $\Exp \colon N_{++} \to S_{++}$ is a diffeomorphism as well.
\end{proof} 
 
 Denote $ {N^i}_{0 j} = \{ (\lambda, t) \in N_6 \mid \tau = i \pi /2, \ \sgn c = j, \ t = 2 \pi / |c| \}, \quad i, j \in \{ +, - \}. $
 
\begin{lmm} \label{ExpNi0j}
Each of the mappings $Exp \colon N_{0 j}^+ \to S_{0 j} , \ \Exp \colon N_{0 j}^- \to S_{0 j} , \ j \in \{ +, - \},$ 
is a diffeomorphism.
\end{lmm}

\begin{proof}
Follows from the parameterization of trajectories~(\ref{ExpN++}) with $p = \pm \pi$.
\end{proof}

Denote $ N_{i 0} = \{ (\lambda, t) \in N_7 \mid \theta = -i \pi/2 \}, \quad i \in \{ +, - \}. $

\begin{lmm} \label{ExpNi0}
The mappings $ \Exp \colon N_{i 0} \to S_{i 0} , \ i \in \{ +, - \},$ are diffeomorphisms.
\end{lmm}

\begin{proof} Follows immediately from the parameterization of extremal trajectories for $\lambda \in C_7$.
\end{proof}

Lemmas~\ref{ExpNij}--\ref{ExpNi0} yield the following optimal synthesis for the terminal manifold $S_6$.

\begin{crllr} \label{S6}
Let $q_1 \in S_6 \backslash \{ q_0 \}$.
\begin{enumerate}
\item If $q_1 \in S_{ij}, \ i, j \in \{ +, - \}$, then the only optimal trajectory is $\Exp(\lambda, t), \ t \in [0, t_1]$,
where $(\lambda, t_1) \in N_{ij} $ is determined by the equations 
$$ x_1 = i \frac{\sin p}{c} , \quad z_1 = \frac{2p - \sin (2 p)}{2 c^2}, \qquad jp \in (0, \pi), \quad jc \in (0, +\infty).$$

\item If $q_1 \in S_{0 j}, \ j\in \{ +, - \} ,$ then there are two optimal trajectories $\Exp({\lambda}_+ , t), \ \Exp({\lambda}_- , t), \ t \in [0, t_1]$, where $({\lambda}_{\pm} , t_1) \in {N^{\pm}_{0 j}}$ is determined by the equations
$$ z_1 = \frac{j \pi}{c^2}, \quad \tau = \pm \frac{\pi}{2}. $$

\item If $q_1 \in S_{i 0}, \ i \in \{ +, -, \}$, then the only optimal trajectory is $\Exp(\lambda, t), \ t \in [0, t_1]$, where $(\lambda, t) \in N_{i 0}$ is determined by the equation $ x_1 = i t_1$.
\end{enumerate}
\end{crllr}

\begin{rmrk} If $\Exp(\lambda, t_1) \in S_6$ for some $(\lambda, t) \in N, $ then ${\varepsilon}^6 (\lambda) = \lambda$.
\end{rmrk}

\begin{proof} It follows from Lemmas~\ref{ExpNij}, \ref{ExpNi0j}, \ref{ExpNi0} that the inclusion $\Exp(\lambda, t_1) \in S_6$ is possible only in 
the following two cases:
\begin{enumerate}
\item $\lambda \in C_6, \ \tau = \theta + ct/2 = \pm \pi/2$, 
\item $\lambda \in C_7, \ \theta = \pm \pi/2$.
\end{enumerate}

In both these cases we have ${\varepsilon}^6 (\lambda) = \lambda$, since ${\varepsilon}^6 \colon (\theta, c, \alpha,) \mapsto (\pi - {\theta}_t , c_t , -\alpha) = (\pi - {\theta}_t, c, \alpha)$ and $\pi - {\theta}_t = \theta.$ 
\end{proof}
\section{Conclusion}
We get a description of the global structure of the exponential mapping in the left-invariant sub-Riemannian problem on the Engel group. It was proved that restriction of this mapping to subdomains in the preimage and image of the exponential mapping cut out by the Maxwell strata corresponding to reflections is a diffeomorphism. Thus we reduced the problem to solving a system of algebraic equations. For any terminal point $q_1 = (x_1, y_1, z_1, v_1)$ with $x_1 \neq 0$ and $z_1 \neq 0$ there exists a unique optimal trajectory. 
Moreover it was proved that the cut time is equal to the first Maxwell time corresponding to reflections.

The cut locus in the sub-Riemannian problem on the Engel group will be described in a forthcoming article. We also plan to study sub-Riemannian spheres and their singularities. Developing software for computation of optimal solutions will allow us to solve the motion planning problem for generic control systems with 4 states and 2 linear inputs via nilpotent approximation (in particular, for the kinematic model of a car with trailer).

\end{document}